What Forms of Student Sensemaking Predict Mathematics Learning?: Student Perspectives on Opportunity, Comfort, and Utility

Sarah E. Castleman, Richard Correnti[1], Christian Schunn[1]

**Abstract**

Instructional practices that require student sensemaking of mathematics are considered critical to learning outcomes, yet strategies for measurement of this instruction at scale have yet to be established. This study aims at validating and exploring the predictive validity of three student survey measures: students' *opportunity* to engage with, *comfort* engaging with, and *utility* of engaging with instructional practices targeting sensemaking in mathematics. Using an initial sample of 1,955 US middle school student responses, we conducted exploratory factor analysis to establish *sensemaking* construct structure, identifying two distinct factors reflecting students' active and passive participation in mathematical sensemaking. We then analyzed data from an additional 3,170 surveys collected over two years in the same context to confirm measurement properties and evaluate the measures' predictive validity. The resulting measures capture students' reported **opportunities** to engage in sensemaking practices (*Sensemaking Opportunity*), their **comfort** with engaging in those practices (*Sensemaking Comfort*), and their perceived **utility** of those practices when tackling new or difficult math problems (*Sensemaking Utility*). While *Sensemaking Comfort* and *Sensemaking Utility* both showed strong predictive validity over two years, *Sensemaking Utility* provided the most reliable estimates of student sensemaking. These findings offer researchers and practitioners a psychometrically sound, scalable toolkit for monitoring sensemaking in mathematics classrooms with implications for both large-scale research on instructional quality and teachers' ongoing efforts to monitor and improve their practice.

[1] The second and third authors contributed equally and are listed in alphabetical order.

## Introduction

Today's professional vision for effective mathematics teaching centers on providing students with opportunities for mathematical sensemaking. Major practice guides (e.g., NCTM, 2014) call on teachers to implement tasks that allow students to explore multiple solution strategies, make connections between representations, engage in productive struggle, and explain their thinking through mathematical discourse—practices identified as essential components of every mathematics lesson (Brahier, Leinwand, & Huinker, 2014). In each of these practices, students are positioned as active participants who work to construct new understanding, often through discourse in which they explain, justify, and co-create meaning with peers (Fitzgerald & Palincsar, 2019; Hiebert et al., 1996). The goal of these practices is conceptual understanding: the capacity to see connections across mathematical ideas and apply them flexibly to new situations, rather than rely on memorized procedures (Schoenfeld, 1991).

Measuring the degree to which these practices are present in classrooms is consequential on multiple fronts: identifying which instructional practices best support student sensemaking, establishing links between sensemaking practices and performance on state assessments, and for providing teachers and researchers with actionable data for iterative improvement.

Yet measuring instructional practices has largely been a time-intensive endeavor. Research in this area has been conducted almost exclusively through qualitative methods (e.g., observations, interviews; Fitzgerald & Palincsar, 2019). These approaches yield rich insights but impose practical constraints: small samples, high resource demands, and limited capacity for real-time use by practitioners. Teacher surveys extend potential reach but come with tradeoffs (e.g., self-report bias).

Student surveys offer promise and have been used to study other areas of instructional practice. Student surveys measuring the frequency of cognitive activation strategies have been linked to self-efficacy, mathematics anxiety, and achievement (Wang, 2023; Zuo, Huang, & Qi, 2024), but have not been designed with mathematical sensemaking as a target construct.

This study extends the measurement literature in three ways. First, researchers have relied on measuring students' visible engagement and mapping this onto students' sensemaking, while acknowledging that students may also be engaging internally or 'covertly' (Chi, 2009) with sensemaking in ways that are invisible to observers. We investigate whether student surveys hold promise for distinguishing between students' active (visible) and passive (invisible) engagement with sensemaking. More specifically, we investigate whether instructional practices related to students' active engagement with sensemaking practices (e.g., explaining strategies) is empirically distinct from passive engagement with sensemaking (e.g., listening to peers' sensemaking). Second, we develop three parallel survey measures of students' reports of the frequency (Opportunity), comfort (Comfort), and perceived helpfulness (Utility) of sensemaking practices, and examine the extent to which these measures offer scalable, psychometrically sound approaches for assessing sensemaking in mathematics classrooms. Third, we compare the predictive validity of all three operationalizations for student mathematics achievement growth, building on prior work linking student-reported instructional practices to outcomes. In doing so, we examine the following research questions:

RQ1: Do items reflecting students' passive participation form a construct distinct from active sensemaking, or integrate into a single dimension?

If two distinct constructs appear:

(1a) Do either *active* or *passive* sensemaking demonstrate predictive validity above and beyond baseline controls?

(1b) Is *active* or *passive* sensemaking a stronger predictor of achievement growth?

RQ2: Can sensemaking be reliably measured in terms of students' perceived *Opportunity* to engage with sensemaking practices, *Comfort* engaging with these practices and *Utility* of engaging with these practices?

RQ3: Do *Utility*, *Opportunity*, and *Comfort* differ in their predictive validity?

We begin with a brief review of the theoretical and empirical literature on mathematical sensemaking and its measurement.

## Theoretical and Empirical Background

### Student Sensemaking as a Goal in Mathematics Instruction

Mathematical sensemaking refers to the process through which students actively construct meaning from mathematical problems, concepts, strategies, and solutions. The student is positioned as a cognitive apprentice, gathering expertise in ways of thinking and doing as an active participant: a maker, a creator of understanding (Lave & Wenger, 1991). Sensemaking is not tied to any specific grade level or mathematical domain but is instead a cross-cutting orientation that students should develop throughout their schooling according to the Common Core (NGA & CCSSO, 2010). Its centrality to high-quality mathematics education is reflected in major practice guides' identification of sensemaking as a hallmark of an excellent mathematics education: students should be engaged in "experiences that promote their ability to make sense of mathematical ideas and reason mathematically" (NCTM, 2014, p. 5).

A defining feature of sensemaking in mathematics is the positioning of the student as actively exploring the subject, seeking to construct, justify and evaluate their own mathematical ideas, as well as those of their peers (Fitzgerald and Palincsar, 2019). Repositioning students as active sense-makers has direct implications for teaching: the role of the teacher shifts from disseminating knowledge into organizing and facilitating a classroom environment in which students can actively make sense of mathematics.

***Active vs. Passive Participation.*** Central to the vision of sensemaking is the role of discussion as a means for students to develop conceptual understanding (English et al., 2025; Jackson et al., 2025; Stein et al., 2008). This pedagogical approach presumes that students build understanding through the interplay of articulating their own thinking, using the verbalized ideas of others to adjust their own understandings, and collectively moving toward a greater shared understanding. This process of advancing the collective knowledge of the group is also known as *building knowledge* (Bereiter & Scardamalia, 2014).

Classroom discussion thus serves a dual function: it creates opportunities for the student who is verbalizing their thinking, but also for the student who is listening — quietly evaluating, confirming, questioning, or revising their own understanding in response. These two forms of participation differ in their visibility. Observable actions such as explaining a strategy or sharing a solution can be documented by teachers and researchers; the cognitive work of the listening student cannot. This asymmetry has shaped the research literature: studies of whole-class discussion in mathematics have focused predominantly on students who are verbally or visibly participating (e.g. Olawoyin, Kribs, & Joswick, 2023), while the sensemaking of listening students remains largely unexamined.

Decades of research indicate that active participation, where students can be seen taking action related to learning (e.g. notetaking, discussing) is preferable to passive participation (e.g. listening to a lecture); instances where students are participating in construction (e.g. self-explanation) or interactive (e.g. productive dialogue) activities with others are better still (Chi, 2009; Chiu & Chi, 2014). Yet relying only on observable engagement makes it impossible to distinguish between a student who is checked-out mentally and one who "may be active internally, mentally self-explaining or making deep connections" (Chiu & Chi, 2014, p. 97). Even research directed specifically at analyzing student listening during mathematics discussion faces this same methodological challenge. After spending a full week observing 100% of one teacher's mathematics instruction, supported by two video-recordings of the teacher and students in every lesson, as well as student- and teacher-interviews, Hintz and Tyson (2015) emphasized their analysis would be based on:

> *if and how students and their teacher responded to one another in one-on-one interactions that took place in front of all students and the teacher. We did not attend to how students and teachers observing the interaction were listening. This limitation is ultimately a methodological challenge when it comes to studying listening – how can a researcher hear what all participants hear? (p. 304)*

Student surveys offer a means of accessing this otherwise invisible dimension of students' experience. By asking students to report on their own engagement as listeners, surveys have the potential to capture forms of sensemaking that observational methods cannot. Building on the self-explanation literature in which students' processing of new information through self-explanations (e.g. making inferences, updating mental models, making connections, expanding on presented information) positions students as active sensemakers (Chi et al., 1989), we hypothesize that this passive participation (i.e. listening to discussion) can constitute a meaningful form of sensemaking. This study uses student surveys to measure and study both students' active and passive participation in mathematical discussion.

## Theoretical Instructional Supports for Student Sensemaking

A substantial body of research explores instructional practices and broader pedagogical approaches that support student mathematical sensemaking. We focus on the practices most closely aligned with positioning students as active sensemakers.

*Solving Problems Independently*. At the most basic level, supporting student sensemaking begins by allowing students to solve problems – "to inquire, to search for solutions, and to resolve incongruities," as opposed to exclusively following routines and procedures that have just been demonstrated (Hiebert et al., 1996, p. 12). Teachers are thus responsible for choosing and presenting tasks that can inspire this sustained inquiry (Hiebert et al., 1996; Stein, Grover & Henningsen, 1996).

*Making & Comparing Multiple Representations.* Students are encouraged to make sense of mathematics (instead of following a prescribed procedure) through exploring multiple solution strategies and invoking multiple representations (e.g., drawing base-ten blocks, vertical subtraction; Hiebert et al., 1996). This approach has been shown to increase students' conceptual understanding: students who were supported and encouraged to first explain the connections between multiple representations (i.e., *sensemaking first*) outperformed students who were first encouraged to notice patterns among specific features of each representations (i.e., *fluency first*); In particular, *sensemaking-first* students showed more gains on fluency and engaged in more conceptual reasoning than *fluency-first* students (Rau, Aleven, Rummel, 2017).

*Providing Justifications.* Justification, or the attempt to justify a mathematical claim using accepted mathematical reasoning, is both a core disciplinary practice as well as a *learning practice*, or a way for students to deepen their understanding of and ability to do mathematics (Cohen & Ball, 2001; Staples, Bartlo & Thanheiser, 2012). Asking students to justify their responses can provide students with access to mathematical discourse and thought, as well as agency in the discipline (Bieda & Staples, 2020).

*Hearing Justifications (i.e., Listening).* There is an assumption in the literature that mathematical justification can occur as part of a social process (i.e. discussion), "But the processes of sense-making and making choices can also be done internally and are not always visible to others" (Lesseig & Lepak, 2022, p. 96). Mercer and Howe (2012) highlight these dual functions of discussion from a Vygotskian lens: the sharing of knowledge through discussion not only contributes to the co-construction of group knowledge (i.e., *collective thinking*), but to *individual thinking* as well: "with children's involvement in joint activities generating new understandings and ways of thinking—not only for them, but *also sometimes for those with whom they are interacting*." (p. 13, emphasis added). Participating in discussions aimed at knowledge production can influence individual thought for both the student sharing as well as the listener. Justification is thus "a learning practice, as it promotes understanding among those engaged in justification—both the individual offering a justification and the audience of that justification—in addition to pursuing a demonstration of the truth of a claim." (Staples, Bartlo & Thanheiser, 2012, p. 449).

*Small and Whole-Group Discussion.* Discussion, whether in small or whole groups, allows students to externalize their thoughts, play out their struggles with learning novel content (English et al., 2025), and to collaboratively build shared mathematical intuition with their peers (Diaz-Rojas, 2027). Discussion grants teachers access to current student thinking and an opportunity to intercede in ways that can redirect, advance or focus student thought (Drageset, 2014), such as by providing explanations necessary to coalesce understanding before the end of class (Gustafsson, 2025). However, these opportunities for teacher intercession are influenced by the teacher's ability to strategically focus on students' thinking (i.e. teacher noticing; Rogers et al., 2025; Seidel et al., 2024), environmental factors (Jazby et al., 2023) and the teacher's ability to pivot from previous plans (Olawoyin, Kribs & Joswick, 2023).

**Typical Measurement Approaches for Students' Sensemaking**

*Teacher Surveys*. Teacher surveys have been used extensively to investigate the instructional practices assumed to support students' mathematical sensemaking at scale (e.g. Banilower et al., 2013; Xia et al., 2025). Findings from these instruments shed light on inequities in students' access to reform-oriented practices (Banilower et al., 2013), and how a teacher's instructional practices can change following professional development (Carney et al., 2015).

Teacher surveys also carry well-documented limitations. For example, teachers' self-reports are frequently misaligned with observational ratings (Kaufman, Stein, & Junker, 2016). Asking teachers to report the frequency of specific practices has been shown to partially address this issue, producing reports more highly correlated with observation-based measures (Mayer, 1999). Weekly teacher logs have also demonstrated reliability for capturing content covered and time spent, though the quality dimension of instruction remains harder to operationalize through self-report alone (Kurz, 2014; Zhang & Schmidt, 2024).

However, teacher surveys measure instruction from the teacher's vantage point (i.e. what teachers report offering) rather than from students' experience of it. This distinction matters beyond issues of self-report bias. The degree to which a given instructional practice creates genuine opportunities for student sensemaking may depend as much on how students engage with it as it does on how the teacher implements it, and teachers may be poorly positioned to report on how well students engage with sensemaking practices. Further, to understand variation within a classroom, teachers are unlikely to be able to report on sensemaking within every individual student or even particular groups of students. Student surveys offer a complementary and underutilized lens, capturing the instructional experiences that students actually perceive and find meaningful.

*Student Surveys.* Student surveys offer a scalable complement to teacher surveys for directly measuring student sensemaking. While student surveys have been underutilized for measuring student sensemaking in particular, they have proved advantageous for investigating related concepts like cognitive activation and meaning-making. Large-scale instruments such as the *Programme for International Student Assessment* (PISA) have captured students' perceptions of the frequency of teachers' use of cognitive activation strategies (e.g. explaining solution approaches) and linked these reports to self-efficacy, mathematics anxiety, and achievement (Wang, 2023; Zuo, Huang, & Qi, 2024).

One line of research spanning over three decades is the work on constructivist-oriented learning environments, most notably through the Constructivist Learning Environment Survey (CLES) and its derivatives (Aldridge et al., 2012; Taylor, Fraser, & Fisher, 1997). Among other aspects of the learning environment, these instruments capture students' perceptions of the frequency of practices associated with collaborative meaning-making (e.g., "I ask other students to explain their ideas) as part of the *student negotiation* scale (Aldridge et al., 2012; Aluri & Fraser, 2019). Iterative versions have been validated across diverse international contexts and used to document instructional change in experimental studies, with some evidence linking constructivist-oriented teaching to pre- and post-test content understanding gains (Ogbuehi & Fraser, 2007). Teachers involved in action research have also found these instruments valuable for reflective practice, using survey results to identify and address gaps in their learning environment (Aldridge et al., 2012). This body of work demonstrates that student surveys can serve as practical, iterative tools for both research and for improving teaching practice.

However, the CLES and related instruments are measures of the broader learning environment rather than of disciplinary sensemaking specifically. Their scope extends well beyond the instructional practices that support students' active construction of mathematical meaning, which limits their utility for identifying which specific sensemaking practices matter most for student learning. A related limitation applies to the cognitive activation literature: while student reports of cognitive activation have been linked to achievement (Wang, 2023; Zuo, Huang, & Qi, 2024), these instruments do not distinguish between students' exposure to sensemaking practices and their engagement with them, a distinction that may matter in understanding more about both student perceptions of instructional practices and how those instructional practices are associated with desirable learning outcomes.

Student surveys have also been used to examine equity dimensions of sensemaking-oriented instruction. Opportunity to learn has been posited to be a key condition connecting teaching to learning, and that this opportunity to learn is particularly present in classrooms where there is explicit attention to mathematical concepts and students have opportunities for productive struggle (Hiebert & Grouws, 2007). Building on this framing, Munter and Haines (2019) used student surveys to investigate whether students from marginalized groups experienced sensemaking-oriented teaching differently from their peers. Students of color were more likely than white students to report racial differences in opportunities to be successful in mathematics. Importantly, among students of color, those in classrooms where teachers initiated cognitively demanding tasks but allowed the cognitive demand of those tasks to decline were significantly more likely to perceive racial discrimination than students of color in other mathematics classrooms. Thus, the presence of sensemaking-oriented tasks is insufficient on its own; students' experience of those tasks matters, and student surveys are better positioned to capture that experience than teacher reports.

Student surveys have been shown to be useful for broad measurement of general learning opportunities for sensemaking. They have even been applied as practical tools for teaching improvement.

In this study, we generated student surveys to measure opportunities for students’ sensemaking specifically in mathematics, and in classrooms where teachers were involved in professional learning to improve their instructional practice. We designed items to explore active as well as passive sensemaking and used various item formats allowing for exploration of the psychometric properties of different variants of students’ sensemaking – i.e., frequency, comfort and utility. While existing surveys have relied primarily on frequency-of-instruction reports, it is an open question whether different question formats provide stronger construct and predictive validity, especially among younger students.

## The present study

The literature reviewed here surfaces three gaps that motivate this study. First, the theoretical distinction between active and passive sensemaking remains empirically unexamined: it is not yet known whether student surveys can capture meaningful information about sensemaking during passive participation (e.g. listening to peers). Second, the field lacks scalable, psychometrically validated instruments for measuring students' sensemaking experiences, which could support both large-scale research and teachers' ongoing improvement work. Third, while student-reported instructional practices have been linked to achievement, it remains unclear which operationalization of sensemaking — frequency, comfort, or utility — best captures the aspects most consequential for learning. This study offers four contributions: empirical evidence on whether active and passive sensemaking constitute distinct dimensions of mathematical engagement; psychometric analysis of scalable survey measures of student sensemaking; initial guidance on which question formats most effectively capture students' sensemaking in middle school mathematics classrooms; and the predictive validity of students' sensemaking on their achievement.

# Methods

## Participants and Context

Participants were middle school students enrolled in three school districts located on the outskirts of a large metropolitan area in the southwestern United States. District A is a mid-sized suburban district with two participating middle schools. Districts B and C are smaller semi-rural and rural districts, respectively. All districts served predominantly Hispanic/Latiné student populations.

This study was conducted as part of a larger networked improvement community (NIC) focused on improving mathematics teaching and learning across the three districts. The NIC provided participating teachers with ongoing professional development and collaborative learning opportunities related to the effective implementation of high-quality curriculum materials, with a special focus on supporting emergent multilingual students.

Teachers participated with the NIC over four school years, from 2021 to 2025. Limited surveys were conducted in the baseline year of 2021. This study examines student survey responses and student achievement growth during three years: (Year 1) the 2022-2023 school year; (Year 2) the 2023-2024 school year; and (Year 3) the 2024-2025 school year. Following recommended practices (Henson & Roberts, 2006), these data were split into two independent samples for the exploratory factor analysis and confirmatory factor analysis phases of scale development.[2]

[2] Due to district-level disruptions (e.g., layoffs, staff turnover, and shifting priorities), limited student enrollment and achievement data were provided in Year 3. However, a large survey dataset was collected during this period ($n$ = 1,955). Accordingly, we used the Year 3 survey data to conduct an exploratory factor analysis (EFA), leveraging the largest available sample to examine the underlying structure of the *Sensemaking* construct. Survey data from Years 1 and 2, when linked achievement data were more complete, were retained for subsequent confirmatory factor analysis (CFA) and analyses of predictive validity.

Across the NIC, 35 middle school mathematics teachers participated in the study, with 15 teachers participating for all three years. Participating teachers taught mainstream mathematics courses (e.g., Grade 6 mathematics; $n$ = 27), special education mathematics courses ($n$ = 3) or Response-to-Intervention (RTI) mathematics courses that students attended in addition to their primary mathematics class ($n$ = 5). On average, participants were highly experienced with an average of more than 14 years in the classroom, reflecting the low teacher turnover in participating districts, particularly as compared to similar districts in the state (TEA, 2025). Participating teacher demographic characteristics are presented in Table 1.

**Table 1.**
*Teacher Demographic Characteristics, Across Study Years*

| | Year 1 | Year 2 | Year 3 | Participated all three years |
|---|---|---|---|---|
| $n$ | 23 | 28 | 26 | 15 |
| Female | 74% | 82% | 81% | 80% |
| Ethnicity | | | | |
| Hispanic or Latiné | 83% | 82% | 81% | 80% |
| Asian | 9% | 7% | 8% | 13% |
| Black | 4% | 4% | 4% | 7% |
| Avg Years Experience | 14.3 | 14.8 | 14.6 | 16.5 |
| Has Master's Degree | 30% | 36% | 35% | 33% |
| Spanish - first language | 39% | 50% | 46% | 40% |

Student participants ($n$ = 5,125) were those who were enrolled in the Grades 6–8 mathematics courses taught by participating teachers, and were given the opportunity to take the student survey at least once. Across all years, over 94% of students identified as Hispanic or Latiné, and more than half reported speaking Spanish at home. Fewer than 5% of students identified as Native American, Black, Asian, or another racial/ethnic group.[3] Table 2 summarizes student demographic characteristics.

## Measures

This study examines the measurement of student-level active sensemaking in mathematics and its association with student-level mathematics achievement growth. The dataset was collected at the student-year level, with each observation reflecting a unique student-teacher instructional relationship within a given academic year. While individual students and teachers may appear in more than one year, each dyad captures a unique instructional relationship at a particular time and under a specific set of classroom learning conditions: a distinct pairing which is therefore treated as an independent relationship.

[3] Due to inconsistencies in district reporting practices regarding the integration of race and ethnicity variables across years, reliable data on students identifying as White were not available for inclusion.

**Table 2.**
*Student Demographic Characteristics, Across Study Years*

| | Year 1 | Year 2 | Year 3 |
|---|---|---|---|
| *n* | 1467 | 1703 | 1955 |
| Grade | | | |
| 6 | 554 | 606 | 647 |
| 7 | 462 | 593 | 610 |
| 8 | 451 | 504 | 698 |
| % Female | 48 | 46 | 49 |
| % Ethnicity | | | |
| Hispanic or Latiné | 93 | 93 | 96 |
| Native American | < 5 | < 5 | < 5 |
| Asian | < 5 | < 5 | < 5 |
| Black | < 5 | < 5 | < 5 |
| % Speaks Spanish at Home | 50 | 49 | 53 |
| % Economic Disadvantage | 65 | 73 | 75 |
| % At-risk | 78 | 64 | 71 |
| % Gifted | 7 | 10 | 12 |
| % Special Education | 10 | 13 | 13 |

*Note.* Descriptive statistics reflect complete observations only; See Appendix A for a summary of missing data and imputation details.

*Mathematics achievement growth.* Students' mathematics achievement growth was measured using the State of Texas Assessments of Academic Readiness (STAAR), a statewide, standards-based and criterion-referenced exam. The test is administered primarily online in a fixed-form (non-adaptive) format. Item types include selected-response (e.g., multiple-choice), technology-enhanced constructed response (e.g., graphing interactions), and open-response numerical items. For analyses, students' mathematics achievement growth was estimated using STAAR percentile ranks at two time points: their percentile rank on the end-of-year STAAR exam ($Math_1$), controlling for their percentile rank on the STAAR exam the previous spring ($Math_0$).

*Student demographics.* Student demographic information was provided by each school district and included each student's mathematics teacher and course enrollment. We constructed indicator (dummy) covariates for each student's demographic information, including gender, ethnicity, home language, participation in the special education and/or gifted programs, economic disadvantage status, and at-risk designation.

*Sensemaking*. The student survey asked students to report on a list of teaching practices in their mathematics class (see Table 3). The list of teaching practices was developed by an interdisciplinary team with expertise in student survey design, mathematics education, bilingual education, and improvement science. Items were designed to reflect instructional practices sponsored through the initiative, incorporate language commonly used by teachers within the initiative, and to be presented in student-friendly language. The teaching practices were independently translated into Spanish by two team members who subsequently collaborated to reconcile discrepancies and produce a final version.

**Table 3.** *Response Distribution of Sensemaking Activities*

| # | Sensemaking Activities | *Opportunity* | | *Comfort* | | *Utility* | |
|---|---|---|---|---|---|---|---|
| | | mean | sd | mean | sd | mean | sd |
| 1 | **Solve** problems on my own | 0.33 | 0.41 | 0.37 | 0.43 | 2.97 | 0.76 |
| 2 | **Create** drawings, tables and graphs | 0.53 | 0.44 | 0.40 | 0.43 | 3.19 | 0.73 |
| 3 | **Compare** drawings, tables, graphs or equations | 0.31 | 0.41 | 0.26 | 0.39 | 3.01 | 0.78 |
| 4 | **Explain** my strategy | 0.72 | 0.39 | 0.41 | 0.44 | 3.03 | 0.79 |
| 5 | **Hear** another student's **strategy** | 0.44 | 0.44 | 0.29 | 0.40 | 3.18 | 0.75 |
| 6 | **Hear** another student **explain** strategy | 0.24 | 0.37 | 0.19 | 0.34 | 3.16 | 0.72 |
| 7 | **Discuss** with a small group or partner | 0.44 | 0.44 | 0.49 | 0.45 | 3.23 | 0.73 |
| 8 | **Share** my strategies with the whole class | 0.32 | 0.41 | 0.19 | 0.35 | 2.64 | 0.83 |
| | Response format | binary | | binary | | Likert (4 point) | |

*Note*. For binary items (*Opportunity, Comfort*), mean reflects the proportion of students responding yes. Table reflects complete cases only.

To examine how middle school students interpret and report on active sensemaking experiences, students were first requested to report on the *perceived Opportunity* and *Comfort with* each of eight teaching practices in their mathematics class. Students were provided with the question stem "My teacher often asks me to…" (*Opportunity*) and "I feel comfortable when I am asked to…" (*Comfort*); and then asked to select as many of the teaching practices (listed in Table 3) as applied. A check-all-that-apply question structure was provided to reduce reading load and completion time. The option "none of the above" was also provided to distinguish between non-response and a lack of frequent exposure to any of the listed practices; surveys with all options blank, including "none of the above," were excluded from analysis. Binary items were created as indicators (1 = selected, 0 = not selected) for each relevant practice for each stem, resulting in two sets of eight binary items named $Opportunity_{\#}$ and $Comfort_{\#}$, respectively.

Next, students were asked about the *utility* of the same teaching practices. They were presented with the question stem, "Which of these activities helped you understand a new or difficult math concept this year?, and asked to respond for each practice using the Likert scale: (1) made me more confused, (2) didn't help at all, (3) helped a little, (4) helped a lot. Eight $Utility_{\#}$ items were thus on a semi-continuous 1-to-4 scale for each of the listed practices.

## Procedures

Student surveys were administered by the mathematics teachers, once in the fall and once in the spring. Students completed the survey in class or as homework. For linking purposes, the survey asked students to provide their student number so that their survey data could be matched with district-provided demographic and performance data. To deal with minor errors in entered numbers, students also indicated their school and teacher.

## Missing Data

To address missingness in survey responses, demographic data, and outcome variables, we used multiple imputation (see Appendix A for details).

## Findings

The findings are organized around three research questions: (1) if *active* and *passive* sensemaking emerge as two distinct constructs and, if so, whether one has stronger predictively validity; (2) whether student survey items capture sensemaking in middle school mathematics as a valid and reliable latent construct, and (3) the predictive validity of these constructs for student mathematics achievement growth.

### Exploring Dimensionality of *Sensemaking: Active* and *Passive* Sensemaking

To address RQ1, we tested whether students' reports of instructional practices in their mathematics classrooms could be explained by an underlying latent construct representing sensemaking in mathematics. We conducted exploratory factor analyses using the most recent survey administration (Year 3), which provided the largest sample and reflected the most stable implementation conditions of the instructional initiative.

Initial exploratory factor analysis investigated whether the items formed a single cohesive latent construct or if students' passive and active participation in sensemaking would separate as distinct constructs. The analyses were conducted on the *Utility* items because these involved a Likert scale to provide the finest item-level sensitivity.

While internal consistency reliability for the eight *Utility* items was strong (McDonald's ω = .83), exploratory factor analysis (EFA) was conducted to assess potential multidimensionality. Eigenvalues and inspection of the scree plot supported a two-factor solution. The first factor had an eigenvalue of 2.8; the second factor had an eigenvalue below 1.0 but corresponded to a visible elbow in the scree plot and explained a fair proportion of variance (.14). Given theoretical expectations that sensemaking would comprise related dimensions, we estimated a two-factor model using a Promax (oblique) rotation.

In the initial solution, all items loaded above .43 on at least one of two factors. However, of the eight items, the item *sharing with a partner or small group* ($Utility_7$) showed weakest loadings. To improve construct coherence and comparability across measures, this item was removed and analysis was repeated.

In the revised solution, all remaining items loaded at least .47 on their primary factor with minimal cross-loadings (see Table 4), revealing a clear distinction between two forms of sensemaking. The first factor included practices in which students actively generated or articulated mathematical ideas, while the second factor included practices in which students engaged with the mathematical reasoning of others.

To evaluate this structure and investigate the independence or covariances of item residuals, a confirmatory factor analysis (CFA) was estimated in Stata using maximum likelihood estimation. The initial two-factor model demonstrated acceptable fit ($\chi^2(13) = 133.19$; CFI = .97, TLI = .95, RMSEA = .07). Modification indices suggested residual covariance between $Utility_2$ (creating representations) and $Utility_3$ (comparing representations). Since both items involve students' engagement with multiple mathematical representations, allowing their residuals to correlate was theoretically justified. Allowing this covariance improved model fit ($\chi^2(12) = 76.69$; $p < .001$; CFI = .98; TLI = .97; RMSEA = .05). Collectively, these findings provide empirical support for modeling students' reports of the helpfulness of instructional practices as a two-dimensional construct reflecting distinct forms of mathematical sensemaking.

**Table 4.**
*Exploratory Factor Analysis of Student Surveys (n=1,955)*

| Variable | Active Sensemaking | Passive Sensemaking |
|---|---|---|
| *Utility₁ :* Solve independently | 0.47 | |
| *Utility₂ :* Create multiple representations | 0.55 | |
| *Utility₃ :* Compare multiple representations | 0.59 | |
| *Utility₄ :* Explain strategy | 0.67 | |
| *Utility₅ :* Hear strategy | | 0.73 |
| *Utility₆ :* Hear explanation | | 0.71 |
| *Utility₈ :* Share whole group | 0.58 | |

*Note.* Factor loadings after oblique Promax rotation. Higher mean scores indicate the practice was more helpful. Loadings below .1 suppressed.

**Construct Validity and Reliability of *Sensemaking Utility***

To examine the stability and measurement properties of this construct, confirmatory factor analyses (CFA) were conducted using Year 1 and Year 2 survey data. Mirroring the exploratory analysis, a two-factor model was estimated for each year, again allowing for covariance between the error terms for creating ($Utility_2$) and comparing ($Utility_3$) multiple representations.

Model fit was strong in both years. To address the limitations of assessing model fit on a multiply imputed dataset, fit indices were first determined on complete cases, which represent actual observed data structure. For Year 1, the model yielded CFI = .99, RMSEA = .03, and SRMR = .03. For Year 2, the model yielded similar specifications: CFI=.98, RMSEA = .05, SRMR = .03. All factor loadings were substantial (≥ .45) and statistically significant ($p < .001$). Model fit was then run on two imputations, which provided similar results and confirmed a well-specified model (see Appendix B). These results indicate good measurement properties and support the stability of the two-factor structure across survey administrations.

To determine whether the two-factor specification provided a better representation of the data than a single-factor alternative, one-model factors were also estimated for each year on complete cases. As shown in Table 5, the two-factor model demonstrated substantially better fit. In Year 1, the two-factor model reduced $\chi^2$ by 229.99 and BIC by 222.88 relative to the single-factor model. Similarly, in Year 2 the two-factor model reduced $\chi^2$ by 324.48 and BIC by 317.11. These improvements indicate that distinguishing between two forms of students' mathematical sensemaking provides a more accurate representation of students' survey responses. Although factor scores during EFA were strongly correlated ($r$ = .84), the substantially better fit of the two-factor model supports their treatment as related but empirically distinct dimensions. To confirm the removal of the item *sharing with a partner or small group* ($Utility_7$), which was removed during EFA, analysis of the two-factor model was repeated in each year with the inclusion of $Utility_7$. These analysis confirmed that in both years, removal of $Utility_7$ resulted in a reduction of $\chi^2$ and improvement of model fit compared to models that included this item.

**Table 5.**
*Comparisons of one- and two-factor models*

| | $\chi^2$ (df) | change in $\chi^2$ | RMSEA | SRMR | CFI | TLI | BIC |
|---|---|---|---|---|---|---|---|
| **Year 1** | | | | | | | |
| 1 factor | 259.41 (13) | | 0.12 | 0.06 | 0.88 | 0.81 | 18542.19 |
| 2 factor | 29.42 (12) | 229.99 | 0.03 | 0.03 | 0.99 | 0.99 | 18319.31 |
| **Year 2** | | | | | | | |
| 1 factor | 391.45 (13) | | 0.14 | 0.06 | 0.86 | 0.78 | 22611.85 |
| 2 factor | 66.97 (12) | 324.48 | 0.05 | 0.03 | 0.98 | 0.97 | 22294.74 |

*Note*. All models allow for covariance between the following two items: creating ($Utility_2$) and comparing ($Utility_3$) multiple representations. Model fit reflects complete cases only.

To evaluate the stability of the factor structure across samples, CFA was conducted independently on each year's dataset. Standardized item loadings were highly consistent across all three years (see Figure 1), with item-level loadings varying by no more than .04 across samples. This stability is notable given that each year represents an independent sample of student-teacher dyads across varying initiative, district, and state-level contextual changes. The replication of both the factor structure, as supported by fit indices, and loading magnitudes across three independent samples provides strong evidence for the robustness of the *Active-* and *Passive- Sensemaking Utility* constructs. Both standardized and unstandardized loadings for Years 1 and 2 complete cases and two imputed cases are reported in Appendix B.

The correlations between the two latent factors were approximately .65 (*p*s < .001) in both years. This indicates a strong relationship between *Active* and *Passive Sensemaking*, consistent with their shared grounding in sensemaking-oriented instruction, while still suggesting that they represent distinct dimensions of students' engagement with mathematical ideas.

## Predictive Validity

### *Base Model Selection*

To examine predictive validity, a baseline HLM model was estimated on multiply imputed data, with students nested in mathematics teachers and survey year included as a dichotomous covariate. The outcome was students' state standardized exam score in mathematics at the end of each academic year ($Math_1$). Covariates were selected based on significance in complete-case preliminary models run separately for each year; prior achievement ($Math_0$), at-risk status, and special education and gifted program participation were retained as controls in the final baseline model due to consistent statistical significance across years. Model structure was evaluated by comparing HLM against a regression with teacher fixed effects. HLM was chosen, as teacher-level variance was meaningful in both years (ICCs of .08 and .14) and the HLM model had a lower BIC.  Analyses were initially conducted independently for each year; a pooled multi-year model with survey year as a fixed effect was ultimately selected to maximize analytic power and precision. Full details on covariate selection and initial structure diagnostics are provided in Appendix C.

**Figure 1.**
Final Model of *Active*- and *Passive-Sensemaking Utility* Across 2 Years (standardized, all loadings p < .001)

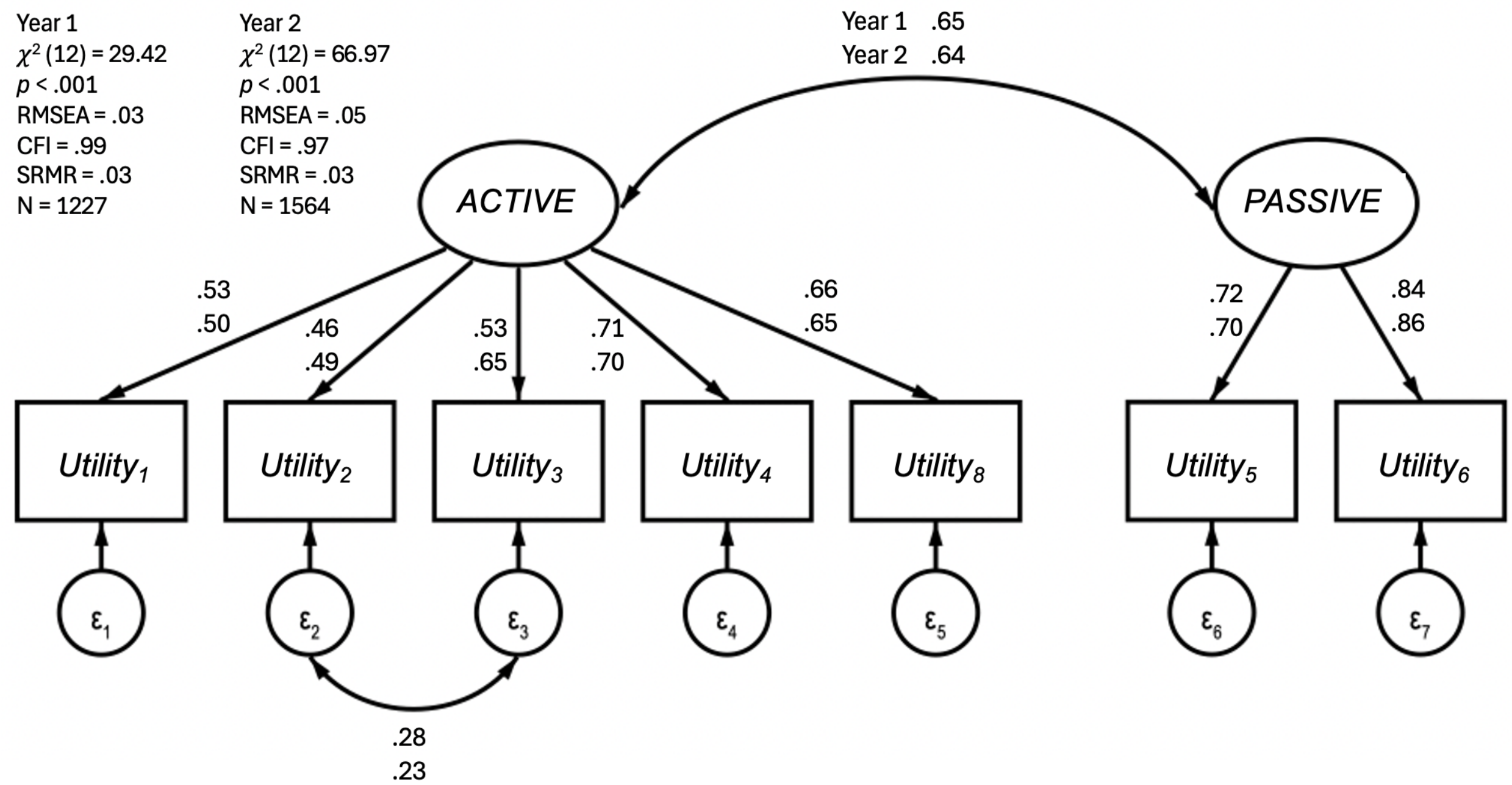

*Note.* CFA was conducted independently on each year's data, using complete cases only. CFA on two imputed datasets is available in Appendix B, with item-level standardized loadings varying by no more than .03 across complete and imputed datasets.

### *Active vs. Passive Sensemaking Utility*

Next, we examined the predictive validity of *Active* sensemaking and *Passive* sensemaking. First, we looked to see if *Active* or *Passive* sensemaking predicted growth above and beyond baseline controls (RQ 1a). When entered in separate models, both constructs significantly predicted student growth, with a larger association for *Active* ($b$ = 3.01, $p$ < .001) than *Passive* sensemaking ($b$ = 2.24, $p$ < .001), as shown in Table 6. This confirmed our hypothesis that *Active* would be a stronger predictor of student mathematics growth than *Passive*, answering RQ 1b. However, when both constructs were entered simultaneously, *Active* retained a significant and similar-magnitude ($b$ = 3.02, $p$ < .001), while the *Passive* covariate had a non-significant relationship ($b$ = -0.02, $p$ = .97) whose coefficient was essentially zero. Taken together, these results support the predictive validity of *Active Sensemaking Utility* and suggest that *Passive* sensemaking does not predict additional student growth beyond its correlation with *Active* sensemaking.

Given these results, we drop *Passive* sensemaking from further analysis and focus subsequent models on *Active* sensemaking alone. Although *Passive* sensemaking was identifiable as a distinct construct and was significantly associated with growth when considered in isolation, once both constructs were entered simultaneously, the coefficient for *Passive* sensemaking was reduced to near zero and lost statistical significance, while *Active* sensemaking remained a robust and stable predictor. This indicates *Passive* sensemaking does not offer additional predictive value

**Table 6.**

*Predictive Validity of Active and Passive Sensemaking Utility (N = 3,170)*

| Model # | 1 | 2 | 3 | 4 |
|---|---|---|---|---|
| | Baseline | Active | Passive | Active & Passive |
| Sensemaking (standardized) | | | | |
| *Active Sensemaking Utility* | | 3.01 *** | | 3.02 *** |
| *Passive Sensemaking Utility* | | | 2.24 *** | -0.02 |
| Achievement | | | | |
| $Math_0$ | 0.58 *** | .56 *** | .57 *** | .56 *** |
| Demographic | | | | |
| At-risk | -6.27 *** | -6.51 *** | -6.46 *** | -6.51 *** |
| Special Education | -5.44 *** | -5.19 *** | -5.07 *** | -5.19 *** |
| Gifted | 9.02 *** | 8.99 *** | 9.08 *** | 8.99 *** |
| Year 1 | 4.55 *** | 4.48 *** | 4.55 *** | 4.48 *** |
| 2 | — | — | — | — |
| Pooled Wald test [1] | | | | |
| over model 1 (baseline) | | 63.81 *** | 34.56 *** | 27.58 *** |
| over model 2 (active) | | | | .00 |
| over model 3 (passive) | | | | 31.29 *** |
| VIF [2] | | | | |
| *Active Sensemaking Utility* | | 1.02 | | 2.55 |
| *Passive Sensemaking Utility* | | | 1.02 | 2.53 |

*Note.* * $p < .05$, ** $p < .01$, *** $p < .001$ — indicates reference category (coefficient constrained to zero) † omitted due to collinearity. [1] Computed with *mi test* in Stata, following Rubin's pooling rules (Li, Raghunathan, Rubin, 1991). [2] VIF computed using OLS on complete cases (*N* = 2,280) as a diagnostic for multicollinearity among fixed effects.

once *Active* sensemaking is accounted for; retaining a covariate with no incremental predictive value would add complexity without analytic benefit. We therefore proceed with Active sensemaking as the focal predictor in the analyses that follow.

## Measurement Reliability of *Sensemaking Utility, Opportunity,* and *Comfort*

To address RQ 2, we next evaluate the construct validity of *Sensemaking* across alternative survey question formats. That is, do items on the five instructional practices also form a coherent construct when measured using students' reports of how frequently they experience these practices (*Opportunity*) and their comfort engaging with them (*Comfort*).

To evaluate this, we estimated confirmatory factor analysis (CFA) models independently on Year 1 and 2 data using the *Opportunity* and *Comfort* items, which were presented as check-all-that-apply questions (i.e., binary response options). Since some students completed multiple surveys within the same mathematics class during a given year, responses were averaged to create a single observation per student per class per year.[4] As a result, *Opportunity* and *Comfort* variables constituted a three-point

[4] Analysis for both Opportunity and Comfort models were repeated at the survey iteration level with dichotomous items, where items similarly showed substantial and statistically significant loadings on the

**Table 7.**

*Means (and SDs) for Survey Items Across Utility, Opportunity, and Comfort Question Stems*

| | *Utility* | | | | *Opportunity* | | *Comfort* | |
|---|---|---|---|---|---|---|---|---|
| | Which of these activities helped you understand a new or difficult math concept this year? | | | | My teacher regularly asks me to… | | I feel comfortable when I am asked to… | |
| | | Year 1 | Year 2 | | Year 1 | Year 2 | Year 1 | Year 2 |
| 1: solve | Solving problems on my own… | 2.96 (.79) | 2.96 (.75) | solve a problem using my own ideas. | .33 (.42) | .33 (41) | .39 (.44) | .36 (.42) |
| 2: make | When I made drawings, tables and graphs, it … | 3.18 (.76) | 3.18 (.71) | make drawings, tables, graphs in solving problems. | .51 (.45) | .55 (.44) | .40 (.44) | .40 (.43) |
| 3: compare | When the teacher had us compare different drawings, tables, graphs or equations, it … | 2.98 (.80) | 3.03 (.77) | compare different drawings, tables, graphs, or equations. | .28 (.40) | .33 (.41) | .26 (.39) | .27 (.39) |
| 4: explain | Explaining why my strategy for solving the problem made sense… | 3.03 (.83) | 3.02 (.77) | explain why my strategy for solving a problem makes sense | .45 (.44) | .43 (.43) | .30 (.41) | .28 (.40) |
| 8: share | Sharing my strategies with the whole class … | 2.65 (.85) | 2.62 (.84) | share my strategies with the whole class. | .34 (.43) | .30 (.40) | .21 (.36) | .18 (.34) |

*Note*. Descriptive statistics based on complete observations; multiple imputation was used in all analyses to address missing data (see Appendix A).

semi-continuous scale (0, 0.5, 1). Descriptive statistics (e.g. mean, SD) of the relevant items across the three question stems are presented in Table 7.

The five-item *Opportunity* scale showed marginally acceptable internal consistency in both Year 1 (McDonald's ω = .68) and Year 2 (ω = .70). As depicted in Figure 2a, confirmatory factor analysis indicated that all items had substantive loadings above .37 ($p < .001$) in Year 1, with good overall model fit ($\chi^2(4) = 21.10$; CFI = .98, TLI = .95, RMSEA = .06). A similar pattern was observed in Year 2, where all items loaded at least .43 ($p < .001$) and model fit remained strong ($\chi^2(4) = 18.42$; CFI = .99, TLI = .97, RMSEA = .05). In order to assess potential impact of missing data, CFA was subsequently performed on two imputed datasets for each year, with similar results (see Appendix B). These results suggest that the *Opportunity* items provide a coherent measure of *Sensemaking*; this *Sensemaking Opportunity* model is presented in Figure 2a.

relevant latent construct. Iteration-level analysis completed in STATA using gsem. Residual covariance could not be estimated for binary indicators.

**Figure 2.**
Final Model of *Sensemaking Opportunity* and *Sensemaking Comfort* Across 2 Years (standardized, all loadings p < .001)

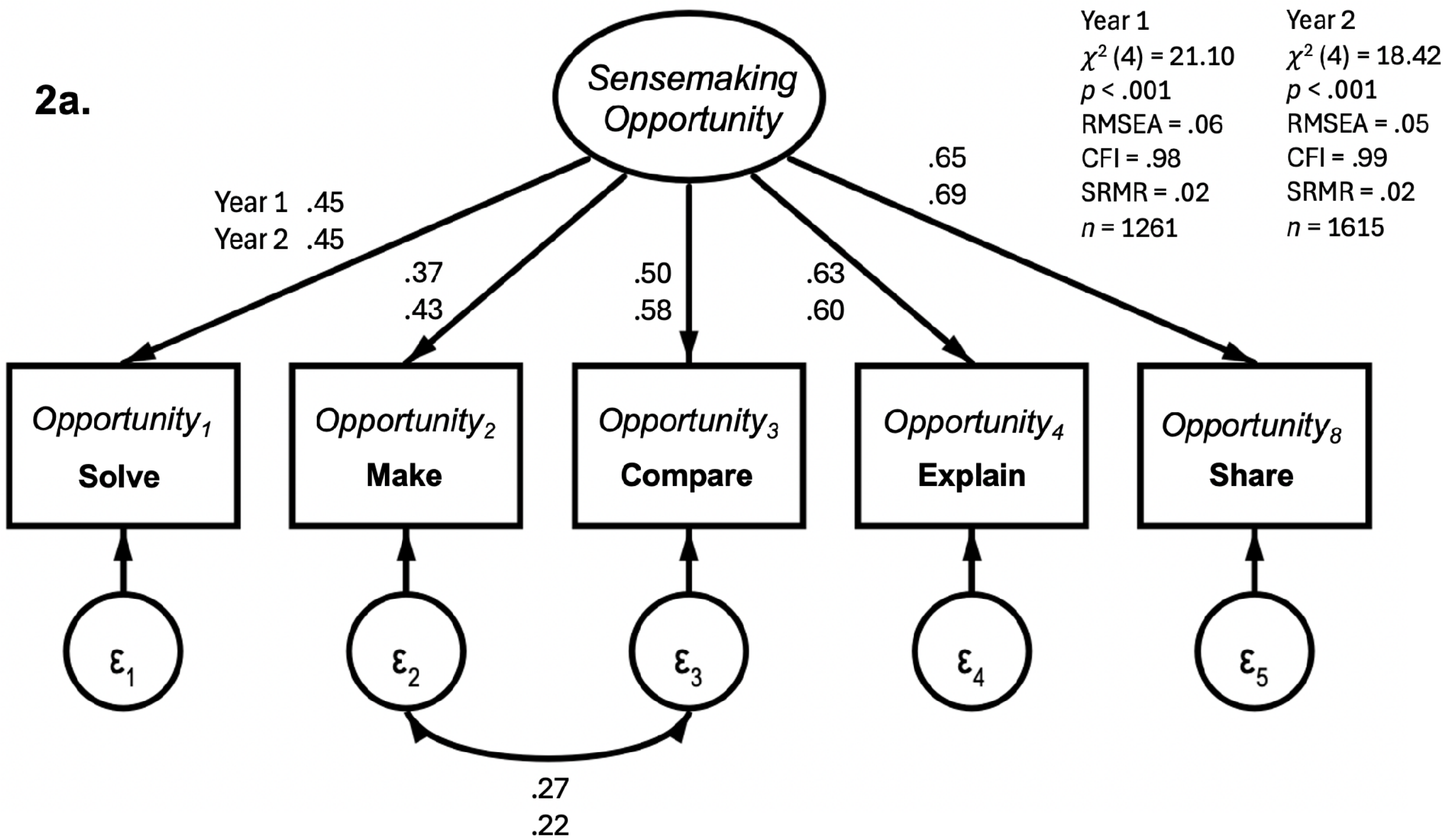

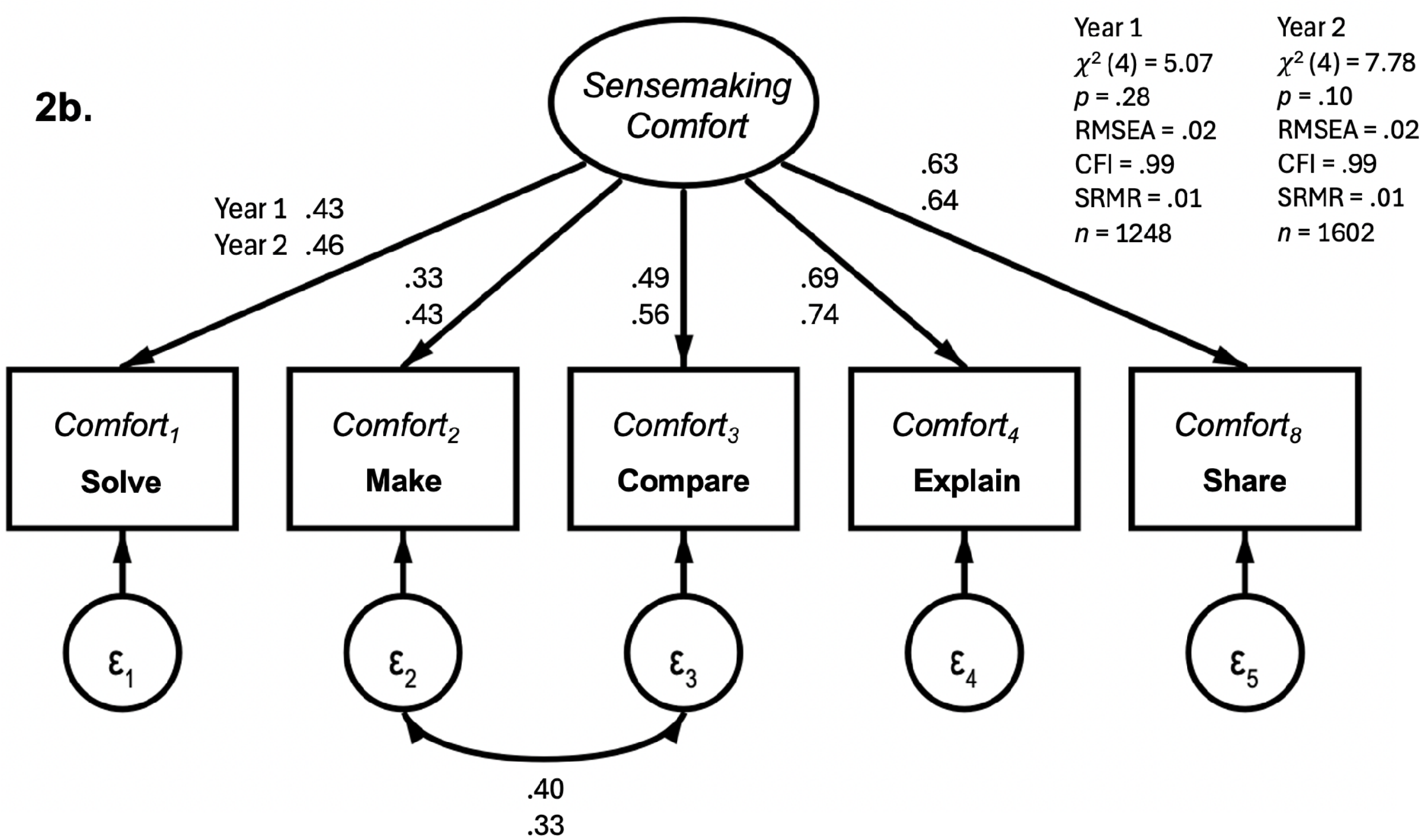

*Note.* CFA was conducted independently on each year's data, using complete cases only. CFA on two imputed datasets is available in Appendix B, with item-level loadings varying by no more than .02 (*Opportunity*) and .04 (*Comfort*) across complete and imputed datasets for Year 1 and Year 2.

The five-item *Comfort* scale showed comparable and acceptable reliability (Year 1 ω = .68; Year 2 ω = .73). In the Year 1 CFA, most items loaded above .43 ($p < .001$), with the exception of *Comfort*$_2$ (i.e., creating multiple representations), which loaded more weakly (.33, $p < .001$). Despite this lower loading, the overall model demonstrated excellent fit ($\chi^2(4) = 5.07$ ; CFI = .99, TLI = .99, RMSEA = .02). Year 2 results showed a similarly good fit: all items had substantive loadings above .43 ($p < .001$) and model fit remained excellent ($\chi^2(4) = 7.78$; CFI = .99, TLI = .99, RMSEA = .02). Replications of model fit on two imputed datasets confirmed model specification (see Appendix B). Together, these findings indicate that the *Comfort* items also capture a coherent latent representation of *Sensemaking Comfort* across years; the final *Sensemaking Comfort* model is depicted in Figure 2b.

Following confirmation of model specification, factor scores were estimated independently for each year across all 20 imputed datasets, and standardized within each year. HLM analyses were conducted within each imputation, and then results were pooled according to Rubin's rules, using 'mi estimate' in Stata 19.5.

Table 7 presents correlations between the four latent sensemaking constructs, achievement variables, and demographic covariates for Years 1 and 2. Among the latent constructs, *Active* and *Passive Sensemaking Utility* were strongly intercorrelated (r = .78/.76), confirming that while they represent distinct constructs, their high shared variance justifies prioritizing *Active Sensemaking Utility* as the sole focal *Sensemaking Utility* construct. *Sensemaking Opportunity* and *Comfort* showed moderate correlations with each other (r = .50/.52), yet each showed only weak correlations with the *Sensemaking Utility* constructs, supporting their treatment as related but distinct operationalizations of sensemaking.

The weak associations between sensemaking constructs and demographic variables reduce concern that observed connections with achievement merely reflect structural confounds. However, the relatively high correlations within some demographic variables suggests that systematic model building may be necessary to avoid multicollinearity issues based within demographic controls. Active Sensemaking Utility showed small negative correlations with Economic Disadvantage (*r* = -.11/-.10) and no significant correlations with at-risk status, home language, or special education participation. This indicates that student reports of *Sensemaking* are not meaningfully correlated with their gender, ethnicity, their language use (*i.e.*, multilingual learners), nor their participation in special education or gifted education. The strong intercorrelations between demographic characteristics, such as those > |.30| and highlighted in green, highlight the possibility for multicollinearity and motivate the need to remove the weakest predictors (*e.g., Male*) during model building to avoid multicollinearity.

The sensemaking constructs showed modest correlations with achievement; *Active Sensemaking Utility* showed the strongest association with end-of-year achievement ($Math_1$: $r = .22/.24$), and a weaker association with prior-year achievement ($Math_0$: $r = .11/.16$), suggesting a closer relationship with current-year growth than with students' pre-existing achievement levels. This pattern would be consistent with the assumption that *Sensemaking Utility* is more reflective of students' engagement with instruction during the year than with students' prior math achievement.

**Table 7.**
*Correlations between Latent Constructs and Potential Covariates (all p<.05)*

| | 1 | 2 | 3 | 4 | Achievement | 5 | 6 | Demographic | 7 | 8 | 9 | 10 | 11 | 12 | 13 | 14 |
|---|---|---|---|---|---|---|---|---|---|---|---|---|---|---|---|---|
| Latent Constructs | | | | | | | | | | | | | | | | |
| 1. *Sensemaking Utility: Active* | — | 0.76 | 0.23 | 0.36 | | 0.24 | 0.16 | | | | | | | | | |
| 2. *Sensemaking Utility: Passive* | 0.78 | — | 0.20 | 0.20 | | 0.20 | | | | | | | | | | |
| 3. *Sensemaking Opportunity* | 0.17 | 0.13 | — | 0.52 | | | | | | | | | | | | |
| 4. *Sensemaking Comfort* | 0.26 | 0.11 | 0.50 | — | | 0.15 | 0.15 | | | | | | | | | |
| Achievement | | | | | | | | | | | | | | | | |
| 5. STAAR1 | 0.22 | 0.14 | | | | — | 0.67 | | | | -0.15 | -0.36 | -0.11 | -0.28 | 0.32 | -0.16 |
| 6. STAAR0 | 0.11 | | | | | 0.61 | — | | | | -0.19 | -0.45 | -0.15 | -0.28 | 0.35 | -0.23 |
| Demographic Variables | | | | | | | | | | | | | | | | |
| 7. Male | | | | | | | | | — | | | | | | | |
| 8. Hispanic | | | | | | | | | 0.10 | — | 0.22 | 0.16 | 0.22 | | | 0.16 |
| 9. Economic Disadvantage | -0.11 | -0.10 | | | | -0.23 | -0.23 | | | 0.14 | — | 0.31 | 0.29 | | -0.12 | 0.25 |
| 10. At-risk | | | | | | -0.34 | -0.33 | | | 0.22 | 0.12 | — | 0.43 | 0.18 | -0.27 | 0.50 |
| 11. Home Language Spanish | | | | | | -0.16 | | | | 0.21 | 0.15 | 0.35 | — | | | 0.69 |
| 12. Special Education | | -0.10 | | | | -0.20 | -0.20 | | 0.11 | | | | | — | -0.11 | |
| 13. Gifted | | | | 0.11 | | 0.31 | 0.26 | | | | -0.13 | -0.17 | | | — | -0.13 |
| 14. English as a Second Language | | | -0.11 | | | -0.15 | | | | 0.17 | 0.15 | 0.34 | 0.72 | | | — |

Note. Lower triangle = Year 1 correlations; upper triangle = Year 2 correlations. Correlations with p ≥ .05 or |r| < .10 are suppressed. Green highlighting indicates correlations greater than |r| < .30.

### *Predictive Validity of Sensemaking Utility, Opportunity and Comfort*

To examine the predictive validity of *Sensemaking Utility, Opportunity* and *Comfort*, we began by testing each construct in a separate model across increasingly specified models: (1) construct only (2) adding prior achievement ($Math_0$) and (3) adding student characteristics and accounting for year fixed effects. The predictive validity of all three constructs proved robust, as each remained statistically significant across every specification tested. Utility, in particular, stood out for producing consistently larger and more stable coefficient estimates than the other constructs (see Appendix D).

To examine which construct retains independent predictive validity in the presence of the others, all three sensemaking constructs were entered simultaneously into the pooled HLM (Table 9). In the combined model, *Sensemaking Utility* retained a significant effect at a magnitude nearly equivalent to its isolated model estimate. *Sensemaking Opportunity* and *Sensemaking Comfort* both became statistically insignificant and shrank substantially in magnitude to near zero. VIF values for all three constructs were low (1.12-1.44), ruling out multicollinearity as an explanation for the loss of significance in *Opportunity* and *Comfort.*[5]

[5] Calculated with regression on complete cases only for VIF diagnostic.

**Table 9.**
*Predicting Mathematics Achievement Growth from Sensemaking Utility, Opportunity, and Comfort (n = 3,170)*

| Model # | 1 | 2 | 3 | 4 | 5 |
|---|---|---|---|---|---|
| | Baseline | *Utility* | *Opportunity* | *Comfort* | *Combined Model* |
| Sensemaking (standardized) | | | | | |
| *Sensemaking Utility* | | 3.01 *** | | | 2.94 *** |
| *Sensemaking Opportunity* | | | 0.73 * | | 0.16 |
| *Sensemaking Comfort* | | | | 1.04 ** | 0.13 |
| Achievement | | | | | |
| *Math* $_0$ | 0.58 *** | .56 *** | .57 *** | .57 *** | .56 *** |
| Demographic | | | | | |
| At-risk | -6.27 *** | -6.51 *** | -6.23 *** | -6.27 *** | -6.49 *** |
| Special Education | -5.44 *** | -5.19 *** | -5.37 *** | -5.43 *** | -5.18 *** |
| Gifted | 9.02 *** | 8.99 *** | 8.98 *** | 8.87 *** | 8.97 *** |
| Year 1 | 4.55 *** | 4.48 *** | 4.57 *** | 4.57 *** | 4.49 *** |
| 2 | — | — | — | — | — |
| Pooled Wald test [1] | | | | | |
| over model 1 (baseline) | | 63.81 *** | 4.04 * | 7.55 ** | 21.94 *** |
| over model 2 (*Utility*) | | | | | 0.21 |
| over model 3 (*Opportunity*) | | | | | 30.02 *** |
| over model 4 (*Comfort*) | | | | | 29.47 *** |
| VIF [2] | | | | | |
| *Sensemaking Utility* | | 1.02 | | | 1.12 |
| *Sensemaking Opportunity* | | | 1.01 | | 1.36 |
| *Sensemaking Comfort* | | | | 1.01 | 1.44 |

* *p* < .05, ** *p* < .01, *** *p* < .001 — indicates reference category (coefficient constrained to zero) [1] Computed with mi test in Stata, following Rubin's pooling rules (Li, Rughunathan, Rubin, 1991). [2] VIF computed using OLS on complete cases (*N* = 2,280) as diagnostic for multicollinearity among fixed effects.

### *Predictive Validity of Sensemaking Constructs vs. Component Items*

To examine whether latent sensemaking constructs offer added predictive value over individual survey items, each predictor (e.g. component item or latent construct) was added individually into the fully specified HLM model controlling for prior achievement and student characteristics. Results are presented in Figure 3; full results are presented in table D1 in the appendix.

For *Utility*, all five component items significantly predicted mathematics achievement growth. One item — Solve independently — did so at a slightly higher magnitude than the *Sensemaking Utility* latent variable (β = 3.34 vs. β = 3.03, p < .001); the remaining items showed smaller standardized coefficients than the latent variable (β = 0.86–2.32). For *Opportunity*, the pattern is notably different: two items failed to reach significance, and no item reached p < .001. The latent variable (β = 0.73, p < .05) consolidated this weak and variable item-level signal into a more stable, if modest,

**Figure 3**

*Predictive Validity of Sensemaking Constructs vs. Component Survey Items (n = 3,170)*

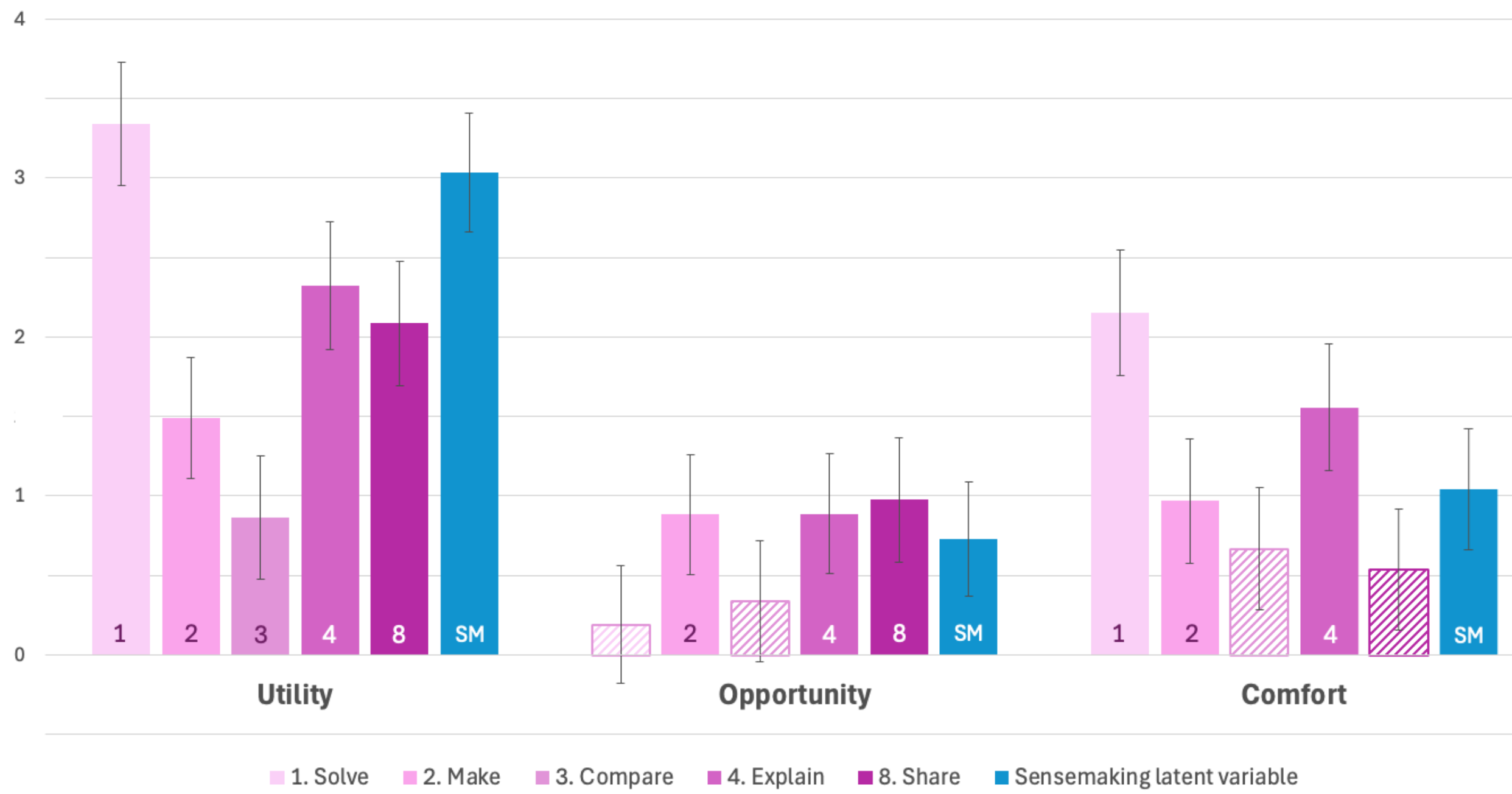

*Note.* Solid bars indicate $p \leq .05$ ; diagonally shaded bars indicate $p > .05$

construct-level estimate. For *Comfort*, three of five items significantly predicted growth. Solve independently again emerged as the strongest individual predictor ($\beta = 2.16$, $p < .001$), performing at approximately twice the magnitude of the *Sensemaking Comfort* latent variable ($\beta = 1.04$, $p < .01$).

Across all three question types, comparing multiple representations showed the weakest predictive performance in both magnitude and significance. *Utility* items collectively showed the strongest and most consistent item-level predictive validity ($\beta$ = 0.86–3.34).

## Discussion

The purpose of this study was to determine if student survey reports of instructional practice could capture student sensemaking as a valid, reliable and predictively meaningful construct.

### The Structure of Sensemaking: Active and Passive Dimensions

A central finding of this study is the empirical differentiation of active and passive sensemaking as distinct, if related, dimensions of students' mathematical experience. The two-factor structure identified through EFA was replicated with minimal variation in factor loadings across three consecutive years. This strengthens confidence that within this context, active and passive sensemaking reflects a stable feature of how students experience sensemaking rather than an artifact of a particular sample or implementation context.

The correlation between the two factors (.64–.65 both years) is itself informative and not surprising. The items share a common instructional basis; both active and passive sensemaking depend on dialogic pedagogical approaches: teachers who implement high-level tasks that welcome student reasoning (Stein, Grover, & Henningsen, 1996) and who facilitate opportunities for students to share their thinking and justify their responses (Staples, Bartlo, & Thanheiser, 2012).

Prior work on measuring cognitive engagement has largely treated passive participation as methodologically inaccessible, a state indistinguishable from disengagement by outside observers (Chi, 2009). A contribution of this study is demonstrating that a passive dimension of sensemaking can be captured through student self-report, consistent with the self-explanation literature's argument that students can be active internal sense-makers even when that participation is not externally visible (Chi et al., 1989; Chi, 2009). In this initial analysis, however, the passive construct did not retain independent predictive validity once active sensemaking was controlled, and its effect not only lost significance but dropped to near zero when the two were entered simultaneously. One interpretation is that passive sensemaking is a consequence of the same instructional conditions that produce active sensemaking rather than an independent contributor to learning: students who report benefiting from listening to peers may simply be those in classrooms rich with sensemaking-oriented instruction, and it is their active engagement in that instruction that drives the achievement relationship. This finding is also consistent with Chi's (2009) framework, which positions passive, active, constructive and interactive as increasingly productive forms of engagement, suggesting that while different students in the same classroom may simultaneously be demonstrating passive listening while another actively shares, it is the latter that drives learning more. The stable passive construct nonetheless holds promise for future work aimed at exploring the cognitive engagement of students who are not providing the externalized evidence of engagement (e.g. discussing, note-taking) visible to outside observers—but who may, nevertheless, be participating in mental processes of mathematical sensemaking.

### *Question Format as a Window into Different Dimensions of Sensemaking*

A second contribution lies in attending to three different question formats (*Utility*, *Opportunity*, and *Comfort*) in building measures of student *sensemaking.* Notably, these provide three related but not equivalent operationalizations of sensemaking. Despite sharing the same item content, the weak inter-construct correlations between *Utility* and both *Opportunity* and *Comfort* indicate these are non-redundant measures. Rather, each measure taps into something different about the student experience.

The field has largely relied on frequency-based question formats for measuring instructional quality, likely because use of frequency reports in teacher surveys has been shown to improve reliability and reduce desirability bias (Mayer, 1999). For student surveys, however, the implications of different question formats have received less attention. Student reports of instruction are not an objective account of what took place in a classroom; their observations of instruction are filtered through students' own engagement with it. A student who is disengaged may report very differently on the frequency of instructional practices than a fully engaged peer in the same classroom, even if the instruction itself is identical.

*Comfort* items add a further layer, reflecting students' affective relationship to instructional practices, shaped by prior experiences, classroom norms, and sense of mathematical identity. *Utility* items ask something different still: a retrospective judgment about what worked—which instructional practices helped the student understand a new or difficult concept.

Both the linguistic difference of the three formats and the distinct aspects of instruction each format taps into may carry particular significance for students. For younger students in particular, these nuances may provide meaningfully different kinds of access to reflecting on and expressing their instructional experiences. This is a dimension of student survey design that warrants further investigation.

*Utility* substantially outperforms *Opportunity* and *Comfort* as a predictor of achievement growth, retaining its effect when all three are controlled simultaneously, while *Opportunity* and *Comfort* shrink to near zero. This suggests *Utility* is capturing something closer to students' actual cognitive engagement with instruction rather than their mere exposure to or comfort with it. This raises questions about the degree to which opportunity-to-learn surveys, which dominate the measurement literature, are capturing what matters most for learning outcomes, at least among middle school students. Where *Opportunity* asks whether the teacher regularly creates conditions for sensemaking, *Utility* asks whether the student takes up those conditions in ways that impact their understanding. Attending to this distinction in student surveys may reveal both differences in instructional quality and in individual student engagement with instruction.

### *Predictive validity*

The study provides a robust empirical link between *Active Sensemaking Utility* and student growth on end-of-year state assessments, even after controlling for prior achievement, at-risk status, and participation in special education and gifted programs. What makes this finding particularly striking is the measure's simplicity: a single question stem with five response options asking students which instructional practices helped them understand a new or difficult mathematics concept. If these findings replicate in other contexts, this constitutes an exceptionally low-burden tool (e.g. quick to administer, straightforward to score) capable of providing researchers and teachers alike with actionable information about which students are engaging with sensemaking-oriented instruction and which may benefit from additional support.

### *Limitations*

The sample is drawn from a networked improvement community, where participating teachers opted-into a program oriented toward high-demand instructional practices; this limits the natural variation in instructional approaches that would exist in other contexts. Data were collected from districts primarily serving Hispanic or Latiné students in a single state; replication in more varied settings, including contexts with greater demographic diversity and different curricular conditions, is needed to confirm the measure's construct and predictive validity.

The student survey measures are not validated with classroom observations, which is the conventional approach to validating teacher reports of instructional quality. While the predictive connection between *Sensemaking Utility* and student achievement is clear, the mechanism is not directly established by this study design. We hypothesize that student reports of sensemaking frequency would correlate with observational measures of sensemaking in the classroom, and that utility and comfort items may reflect students' actual engagement with sensemaking practices, but these hypotheses require further examination.

The passive sensemaking factor comprises only two items, which may limit its predictive validity. Future work should examine whether a broader set of passive participation items, capturing more varied forms of (invisible) cognitive processes occurring during dialogic teaching, could produce stronger and more differentiated findings.

### *Implications for Research and Practice*

These findings carry implications across several dimensions of research and practice.

*As a large-scale research instrument.* Research on teaching practices that support sensemaking has typically relied on time-intensive observational methods that necessitate focus on one or a small number of teachers (Fitzgerald & Palincsar, 2019). The validated *Sensemaking Utility* construct offers a single-item format that can easily be deployed across large and diverse samples to examine the distribution of sensemaking-oriented instruction, its relationship to student outcomes, and its variation across demographic groups and instructional contexts.

*As a practical tool for iterative improvement.* For any context where teachers or teacher-researchers may be reflecting on their teaching (e.g. teacher-coach pairs, grade-level teams), a representation of instruction is needed to anchor productive discussion (Horn, 2010; Horn, Kane, & Wilson, 2015). The measures studied in this paper offer student reports of instruction which are low-investment to collect, interpretable without specialized training, and capable of revealing meaningful differences in how students are experiencing instruction. This combination also positions these measures well for action research; providing teachers with opportunities to conduct action research in their own classroom (e.g., testing replication of empirical findings) can increase teacher buy-in and potentially be a road to instructional change (Aldridge et al., 2012). The construct's predictive relationship with state assessment performance also provides a practical basis for buy-in from school and district leaders invested in instructional changes that translate to measurable student growth.

*As a means of centering student voice.* Perhaps most distinctively, student surveys of sensemaking shift the lens of formative assessment from evaluation of student performance to students' own reports of their instructional experience. Students are positioned as the experts on instructional utility, the informants reporting on their own relationship to the sensemaking practices of their classroom. This reframing has both practical and equity implications: it has the potential to remove mathematics anxiety from the survey experience while still generating predictive signal for achievement growth, and it situates students as agentic contributors to the improvement of instruction rather than as passive objects of assessment.

## Conclusion

This study set out to examine whether student survey reports of instructional practices could capture mathematical sensemaking as a valid, reliable, and predictively meaningful construct. The answer was yes; the *Utility*, *Opportunity* and *Comfort Sensemaking* measures demonstrated strong validity over three samples, with *Utility* ultimately providing the most robust links to student achievement. *Utility*, derived from students' retrospective judgments of which practices helped them understand difficult mathematics, emerged as a robust, stable, and practically consequential measure. *Utility* predicted mathematics achievement growth above and beyond prior achievement and student characteristics, outperformed *Opportunity* and *Comfort* formats both in isolation and in combination, and did so through an instrument brief enough for routine classroom use.

Three contributions stand out. First, the study provides the field with a validity argument for survey items measuring student sensemaking. Second, it advances the theoretical understanding of sensemaking's structure by demonstrating empirically that active and passive forms of engagement constitute distinct dimensions, offering a foundation for future work on how students' internal engagement during discussion can be captured and understood. Third, it offers initial and practically consequential evidence that how we ask students about their sensemaking matters: students' judgments of the *utility* of instructional practices carry predictive signal that *frequency (i.e., opportunity to learn)* and *comfort* formats do not, suggesting that the field's reliance on opportunity-to-learn surveys may underestimate what student surveys can reveal about the relationship between instruction and learning.

## Appendix A

## Multiple Imputation of Pooled Dataset

### Missing Data

*Survey completion.* Survey completion rates were high in both years, with 87% of students in participating classes completing the survey in Year 1 and 95% in Year 2.[6] However, comparisons between survey completers and non-completers revealed systematic differences in both years, including non-completers consistently showing lower $Math_1$ scores and higher at-risk rates (both $p$<.05). To address these observed patterns, multiple imputation was performed, and the variables associated with missingness were incorporated into the imputation model.

*Missing demographic or performance data.* Of the 3,170 participating students, 447 (14%) were missing spring STAAR scores ($Math_1$), and 438 (14%) were missing prior-year STAAR scores ($Math_0$). Less than 9% of students were missing demographic information on any particular item. To reduce potential bias associated with missing data and to retain analytic power, we employed multiple imputation using chained equations (MICE) in Stata v18.5. Imputation models were estimated separately by year and incorporated all variables used in the analytic models as well as additional demographic characteristics.

### Multiple Imputation

Multiple imputation was performed on the pooled Year 1 and Year 2 dataset to address survey missingness. Prior to imputation, observations missing all survey items and mathematics achievement data were dropped ( $N$ = 35). Imputation was conducted using chained equations (MICE) in Stata v18.5 with 20 imputations. Auxiliary variables (i.e. those with complete observations included to inform the imputation model but not themselves imputed) included survey year, school, grade, and mathematics teacher ID. Additional demographic and achievement variables were included in the imputation model based on statistical significance in preliminary complete-case models (see Appendix C); these are listed in Table A1. Variables used in construct validity (RQ2) or predictive validity (RQ3) analyses are indicated in bold and shaded in Table A1.

**Table A1.**
*Multiple Imputation of Pooled Dataset ( N = 3170 )*

| | Complete observations | Imputed observations | Imputed |
|---|---|---|---|
| Auxiliary Variables | | | |
| **Survey Year** | **3170** | **0** | **0** |
| School | 3170 | 0 | 0 |
| Grade | 3170 | 0 | 0 |
| **Mathematics Teacher** | **3170** | **0** | **0** |

[6] Missingness analyses were conducted with Year 1 and 2 data, as complete student enrollment data were unavailable for a portion of Year 3 participants. However, the ratio of students to teachers indicates the missingness rates were low in Year 3.

| | | | |
|---|---|---|---|
| Imputed Variables | | | |
| Demographic Variables | | | |
| Gender | 3123 | 47 | 1 |
| Hispanic | 3117 | 53 | 2 |
| **At-risk** | **2982** | **188** | **6** |
| Economic Disadvantage | 3020 | 150 | 5 |
| Home Language Spanish | 3023 | 147 | 5 |
| School Program Participation | | | |
| Special Education | 3080 | 90 | 3 |
| **Gifted** | **2906** | **264** | **8** |
| English Second Language | 2977 | 193 | 6 |
| Mathematics Achievement | | | |
| ***Math***$_{1}$ | **2723** | **447** | **14** |
| *Math*$_{1}$ Grade Level | 2727 | 443 | 14 |
| ***Math***$_{0}$ | **2732** | **438** | **14** |
| *Math*$_{-1}$ | 875 | 2295 | 72 |
| Survey Items | | | |
| ***Utility***$_{1}$ | **2845** | **325** | **10** |
| ***Utility***$_{2}$ | **2844** | **326** | **10** |
| ***Utility***$_{3}$ | **2842** | **328** | **10** |
| ***Utility***$_{4}$ | **2846** | **324** | **10** |
| ***Utility***$_{5}$ | **2836** | **334** | **11** |
| ***Utility***$_{6}$ | **2828** | **342** | **11** |
| ***Utility***$_{8}$ | **2826** | **344** | **11** |
| ***Opportunity***$_{1}$ | **2876** | **294** | **9** |
| ***Opportunity***$_{2}$ | **2876** | **294** | **9** |
| ***Opportunity***$_{3}$ | **2876** | **294** | **9** |
| ***Opportunity***$_{4}$ | **2876** | **294** | **9** |
| ***Opportunity***$_{8}$ | **2876** | **294** | **9** |
| ***Comfort***$_{1}$ | **2850** | **320** | **10** |
| ***Comfort***$_{2}$ | **2850** | **320** | **10** |
| ***Comfort***$_{3}$ | **2850** | **320** | **10** |
| ***Comfort***$_{4}$ | **2850** | **320** | **10** |
| ***Comfort***$_{8}$ | **2850** | **320** | **10** |

*Note.* **Bold and highlight indicates items included in final analysis**

## Appendix B

### CFA Model Item Loadings and Model Fit of *Sensemaking Utility*

**Table B1 .**

| | | Year 23 | | | | | | Year 24 | | | | | |
|---|---|---|---|---|---|---|---|---|---|---|---|---|---|
| | | complete cases | | Imputed dataset #1 | | Imputed dataset #20 | | complete cases | | Imputed dataset #1 | | Imputed dataset #20 | |
| Item | Latent construct | β | B | β | B | β | B | β | B | β | B | β | B |
| *Utility*$_1$ | ← *ACTIVE* | 0.53 | — | 0.53 | — | 0.52 | — | 0.50 | — | 0.50 | — | 0.50 | — |
| | | (.03) | | (.02) | | (.02) | | (.02) | | (.02) | | (.02) | |
| *Utility*$_2$ | ← *ACTIVE* | 0.46 | 0.84 | 0.45 | 0.83 | 0.46 | 0.85 | 0.49 | 0.91 | 0.50 | 0.96 | 0.49 | 0.93 |
| | | (.03) | (.07) | (.03) | (.07) | (.03) | (.07) | (.02) | (.07) | (.02) | (.07) | (.02) | (.07) |
| *Utility*$_3$ | ← *ACTIVE* | 0.53 | 1.03 | 0.52 | 1.01 | 0.52 | 1.02 | 0.65 | 1.33 | 0.64 | 1.33 | 0.64 | 1.33 |
| | | (.03) | (.08) | (.02) | (.07) | (.02) | (.07) | (.02) | (.09) | (.02) | (.08) | (.02) | (.08) |
| *Utility*$_4$ | ← *ACTIVE* | 0.71 | 1.42 | 0.70 | 1.41 | 0.70 | 1.41 | 0.70 | 1.43 | 0.70 | 1.45 | 0.69 | 1.42 |
| | | (.02) | (.10) | (.02) | (.09) | (.02) | (.09) | (.02) | (.09) | (.02) | (.09) | (.02) | (.08) |
| *Utility*$_8$ | ← *ACTIVE* | 0.66 | 1.35 | 0.63 | 1.28 | 0.64 | 1.30 | 0.65 | 1.43 | 0.64 | 1.43 | 0.63 | 1.41 |
| | | (.02) | (.09) | (.02) | (.08) | (.02) | (.08) | (.02) | (.09) | (.02) | (.09) | (.02) | (.09) |
| *Utility*$_5$ | ← *PASSIVE* | 0.72 | — | 0.71 | — | 0.71 | — | 0.70 | — | 0.70 | — | 0.70 | — |
| | | (.02) | | (.02) | | (.02) | | (.02) | | (.02) | | (.02) | |
| *Utility*$_6$ | ← *PASSIVE* | 0.84 | 1.13 | 0.84 | 1.16 | 0.83 | 1.13 | 0.86 | 1.16 | 0.86 | 1.17 | 0.86 | 1.17 |
| | | (.02) | (.07) | (.02) | (.07) | (.02) | (.06) | (.02) | (.06) | (.02) | (.06) | (.02) | (.06) |
| cov 2 *3 | | 0.28 | 0.13 | 0.26 | 0.11 | 0.27 | 0.12 | 0.23 | 0.08 | 0.22 | 0.08 | 0.22 | 0.08 |
| | | (.03) | (.02) | (.03) | (.01) | (.03) | (.01) | (.03) | (.01) | (.03) | (.01) | (.03) | (.01) |
| cov ACTIVE * PASSIVE | | 0.65 | 0.14 | 0.65 | 0.13 | 0.66 | 0.14 | 0.64 | 0.13 | 0.64 | 0.12 | 0.64 | 0.12 |
| | | (.03) | (.01) | (.03) | (.01) | (.03) | (.01) | (.02) | (.01) | (.02) | (.01) | (.02) | (.01) |
| chi2 | | 29.42 | | 39.32 | | 35.77 | | 66.97 | | 61.99 | | 66.01 | |
| RMSEA | | 0.03 | | 0.04 | | 0.04 | | 0.05 | | 0.05 | | 0.05 | |
| BIC | | 18319.31 | | 21211.84 | | 21251.68 | | 22294.74 | | 24005.68 | | 23944.16 | |
| CFI | | 0.99 | | 0.99 | | 0.99 | | 0.98 | | 0.98 | | 0.98 | |
| TLI | | 0.99 | | 0.98 | | 0.98 | | 0.97 | | 0.97 | | 0.97 | |
| SRMR | | 0.03 | | 0.03 | | 0.03 | | 0.03 | | 0.03 | | 0.03 | |
| n | | 1227 | | 1467 | | 1467 | | 1564 | | 1703 | | 1703 | |

all $p < 0.001$. B = unstandardized loading; β = standardized loading. For unstandardized loadings, the first indicator for each latent factor was constrained to 1 for model identification. Standard errors are provided in parentheses.

**CFA Model Item Loadings and Model Fit of *Sensemaking Opportunity***

**Table B2** .

| | Year 23 | | | | | | Year 24 | | | | | |
|---|---|---|---|---|---|---|---|---|---|---|---|---|
| | complete cases | | Imputed dataset #1 | | Imputed dataset #20 | | complete cases | | Imputed dataset #1 | | Imputed dataset #20 | |
| Item | β | B | β | B | β | B | β | B | β | B | β | B |
| *Opportunity*$_1$ | 0.45 | — | 0.47 | — | 0.47 | — | 0.45 | — | 0.45 | — | 0.45 | — |
| | (.03) | | (.03) | | (.03) | | (.03) | | (.02) | | (.02) | |
| *Opportunity*$_2$ | 0.37 | 0.86 | 0.36 | 0.80 | 0.36 | 0.80 | 0.43 | 1.04 | 0.43 | 1.01 | 0.43 | 1.00 |
| | (.03) | (.10) | (.03) | (.09) | (.03) | (.09) | (.02) | (.10) | (.03) | (.09) | (.03) | (.09) |
| *Opportunity*$_3$ | 0.50 | 1.07 | 0.52 | 1.07 | 0.52 | 1.07 | 0.58 | 1.31 | 0.58 | 1.30 | 0.58 | 1.30 |
| | (.03) | (.10) | (.03) | (.09) | (.03) | (.09) | (.02) | (.10) | (.02) | (.10) | (.02) | (.10) |
| *Opportunity*$_4$ | 0.63 | 1.47 | 0.62 | 1.37 | 0.62 | 1.37 | 0.60 | 1.42 | 0.60 | 1.40 | 0.60 | 1.39 |
| | (.03) | (.13) | (.03) | (.11) | (.03) | (.11) | (.02) | (.11) | (.02) | (.11) | (.02) | (.10) |
| *Opportunity*$_8$ | 0.65 | 1.46 | 0.67 | 1.43 | 0.66 | 1.43 | 0.69 | 1.50 | 0.69 | 1.49 | 0.69 | 1.49 |
| | (.03) | (.13) | (.03) | (.11) | (.03) | (.11) | (.02) | (.11) | (.02) | (.11) | (.02) | (.11) |
| cov 2 *3 | 0.27 | 0.04 | 0.25 | 0.03 | 0.25 | 0.03 | 0.22 | 0.03 | 0.22 | 0.03 | 0.22 | 0.03 |
| | (.03) | (.005) | (.03) | (.004) | (.03) | (.004) | (.03) | (.004) | (.03) | (.004) | (.03) | (.004) |
| chi2 | 21.10 | | 29.46 | | 29.78 | | 18.42 | | 19.42 | | 19.85 | |
| RMSEA | 0.06 | | 0.07 | | 0.07 | | 0.05 | | 0.05 | | 0.05 | |
| BIC | 6419.55 | | 6465.83 | | 6468.13 | | 7732.86 | | 7744.94 | | 7748.69 | |
| CFI | 0.98 | | 0.98 | | 0.98 | | 0.99 | | 0.99 | | 0.99 | |
| TLI | 0.95 | | 0.94 | | 0.94 | | 0.97 | | 0.97 | | 0.97 | |
| SRMR | 0.02 | | 0.03 | | 0.03 | | 0.02 | | 0.02 | | 0.02 | |
| n | 1261 | | 1467 | | 1467 | | 1615 | | 1703 | | 1703 | |

all $p < 0.001$. B = unstandardized loading; β = standardized loading. For unstandardized loadings, the first indicator for each latent factor was constrained to 1 for model identification. Standard errors are provided in parentheses.

**CFA Model Item Loadings and Model Fit of *Sensemaking Comfort***

**Table B3 .**

| | Year 23 | | | | | | Year 24 | | | | | |
|---|---|---|---|---|---|---|---|---|---|---|---|---|
| | complete cases | | Imputed dataset #1 | | Imputed dataset #20 | | complete cases | | Imputed dataset #1 | | Imputed dataset #20 | |
| Item | β | B | β | B | β | B | β | B | β | B | β | B |
| *Comfort1* | 0.43 | — | 0.44 | — | 0.44 | — | 0.46 | — | 0.47 | — | 0.47 | — |
| | (.03) | | (.03) | | (.03) | | (.02) | | (.02) | | (.02) | |
| *Comfort2* | 0.33 | 0.78 | 0.34 | 0.78 | 0.35 | 0.79 | 0.43 | 0.94 | 0.43 | 0.93 | 0.43 | 0.93 |
| | (.03) | (.10) | (.03) | (.09) | (.03) | (.09) | (.03) | (.08) | (.02) | (.08) | (.02) | (.08) |
| *Comfort3* | 0.49 | 1.02 | 0.53 | 1.10 | 0.53 | 1.10 | 0.56 | 1.12 | 0.58 | 1.14 | 0.58 | 1.13 |
| | (.03) | (.10) | (.03) | (.09) | (.03) | (.09) | (.02) | (.08) | (.02) | (.08) | (.02) | (.08) |
| *Comfort4* | 0.69 | 1.51 | 0.69 | 1.50 | 0.69 | 1.50 | 0.74 | 1.51 | 0.74 | 1.50 | 0.74 | 1.49 |
| | (.03) | (.14) | (.02) | (.12) | (.02) | (.12) | (.02) | (.11) | (.02) | (.10) | (.02) | (.10) |
| *Comfort8* | 0.63 | 1.21 | 0.66 | 1.30 | 0.66 | 1.30 | 0.64 | 1.11 | 0.65 | 1.14 | 0.65 | 1.14 |
| | (.03) | (.11) | (.02) | (.10) | (.02) | (.10) | (.02) | (.08) | (.02) | (.07) | (.02) | (.07) |
| cov 2 *3 | 0.40 | 0.06 | 0.40 | 0.05 | 0.40 | 0.05 | 0.33 | 0.04 | 0.33 | 0.04 | 0.33 | 0.04 |
| | (.03) | (.004) | (.03) | (.003) | (.02) | (.003) | (.03) | (.003) | (.03) | (.004) | (.02) | (.004) |
| chi2 | 5.07 | | 4.19 | | 4.00 | | 7.78 | | 8.28 | | 8.12 | |
| RMSEA | 0.02 | | 0.01 | | 0.01 | | 0.02 | | 0.03 | | 0.03 | |
| BIC | 5594.08 | | 5569.85 | | 5578.87 | | 6395.86 | | 6379.16 | | 6374.37 | |
| CFI | 0.99 | | 0.99 | | 0.99 | | 0.99 | | 0.99 | | 0.99 | |
| TLI | 0.99 | | 0.99 | | 0.99 | | 0.99 | | 0.99 | | 0.99 | |
| SRMR | 0.01 | | 0.01 | | 0.01 | | 0.01 | | 0.01 | | 0.01 | |
| n | 1248 | | 1467 | | 1467 | | 1602 | | 1703 | | 1703 | |

all $p < 0.001$. B = unstandardized loading; β = standardized loading. For unstandardized loadings, the first indicator for each latent factor was constrained to 1 for model identification. Standard errors are provided in parentheses.

## Appendix C

### Preliminary Complete-Case Model Diagnostics: Model Structure

**Table C1.**

| | Year 23 | | Year 24 | | Pooled model | |
|---|---|---|---|---|---|---|
| Enrollment | Regression | HLM | Regression | HLM | Regression | HLM |
| Teacher-level effects | Teacher dummies | Random intercepts | Teacher dummies | Random intercepts | Teacher dummies | Random intercepts |
| Grade | | | | | | |
| 6 (reference) | — | — | — | — | — | — |
| 7 | 1.65 | 2.83 | -13.71 | -13.39 * | 1.88 | 1.57 |
| 8 | 1.85 | 1.78 | -5.81 | -8.60 | 3.99 | 3.05 |
| Achievement | | | | | | |
| *Math*$_0$ | 0.62 *** | 0.61 *** | 0.54 *** | 0.54 *** | .57 *** | .58 *** |
| *Math*$_1$ Grade | | | | | | |
| 6 (reference) | — | — | — | — | — | — |
| 7 | 1.08 | 2.16 | 9.12 * | 3.80 | -0.09 | -0.12 |
| 8 | 10.88 * | 9.00 * | † | † | 3.65 | 3.44 |
| Student characteristics | | | | | | |
| Male | 0.18 | -0.15 | 1.92 | 1.93 | 0.97 | 0.90 |
| Race / Ethnicity | | | | | | |
| AmInd | -2.64 | -2.69 | 4.04 | 4.20 | -0.19 | -0.19 |
| Asian | 5.62 | 5.64 | -2.79 | -3.12 | 1.02 | 0.99 |
| Black | - 8.08 | -7.53 | 0.23 | 0.40 | -3.65 | -3.61 |
| PacIsl | 9.21 | 9.26 | 14.96 | 14.70 | 9.33 | 9.02 |
| Hisp | -1.07 | -0.72 | -1.97 | -2.85 | -1.56 | -1.75 |
| Migrant | -3.61 | -3.64 | 4.31 | 4.41 | 1.05 | 1.17 |
| At-risk | -6.05 ** | -6.10 ** | -5.00 * | -4.65 * | -5.21 *** | -4.83 *** |
| Economic Disadvantage | -5.27 *** | -5.56 *** | 0.50 | 0.61 | -1.35 | -1.43 |
| Home Language Spanish | 3.15 | 2.54 | -0.84 | -1.48 | 2.23 | 1.64 |
| English Second Language | -0.40 | -0.31 | 1.55 | 1.55 | -0.22 | -0.28 |
| Emergent Bilingual | -3.98 | -4.42 | 0.04 | 0.81 | -2.64 | -2.11 |
| Special Education | -2.98 | -3.80 | -9.29 *** | -9.87 *** | -6.02 *** | -6.50 *** |
| Gifted | 9.00 *** | 9.42 *** | 3.10 | 3.37 | 1.67 | 1.86 |
| Year | | | | | | |
| 1 | | | | | 1.73 | 2.03 |
| 2 | | | | | — | — |
| Gifted x Year 1 | | | | | 8.08 * | 8.06 * |
| Gifted x Year 2 | | | | | † | † |
| n | 905 | 905 | 784 | 784 | 1689 | 1689 |
| ICC | | .083 | | .144 | | .126 |
| r2 | 0.51 | | 0.56 | | 0.51 | |
| BIC | 8011.71 | 7948.59 | 6906.77 | 6846.63 | 14848.71 | 14749.23 |

*Note.* * $p < .05$, ** $p < .01$, *** $p < .001$ † omitted due to collinearity. Race/ethnicity, emergent bilingual and migrant status were excluded from the multiple imputation model due to consistently insignificant effects across preliminary models; individual race/ethnicity indicators (e.g. Black) were additionally subject to perfect prediction errors during imputation.

**Preliminary Imputed Datasets ($n$ = 20) Model Diagnostics:
Covariate selection ( $n$ = 3,170)**

**Table C1.**

| Model | 1 | 2 | 3 |
|---|---|---|---|
| | All var included in multiple imputation | Final model | |
| Grade | | | |
| 6 (reference) | — | | |
| 7 | -0.07 | | |
| 8 | 1.69 | | |
| Achievement | | | |
| *Math* $_0$ | 0.57 *** | 0.58 *** | |
| *Math* $_1$ Grade | | | |
| 6 (reference) | — | | |
| 7 | 1.03 | | |
| 8 | 2.92 | | |
| Student characteristics | | | |
| Male | 0.98 | | |
| Hispanic | -1.20 | | |
| At-risk | -5.48 *** | -6.27 *** | |
| Economic Disadvantage | -1.78 | | |
| Home Language Spanish | 0.07 | | |
| English Second Language | -1.27 | | |
| Special Education | -5.83 *** | -5.44 *** | |
| Year | | | |
| 1 | 3.52 *** | 4.55 *** | |
| 2 | | | |
| Gifted x Year 1 | 4.99 | | |
| Gifted x Year 2 | † | | |

*Note.* * $p$ < .05, ** $p$ < .01, *** $p$ < .001 † omitted due to collinearity.

Economic Disadvantage and Gifted x Year interaction were initially tested due to significance in preliminary complete case analysis, but ultimately dropped due to insignificance in the HLM model on imputed data.

## Appendix D

***Predictive Validity of Sensemaking Across Model Specifications***

To assess robustness of the three *Sensemaking* constructs, each latent construct was tested across three increasingly specified models: (1) construct only, with year as a covariate; (2) adding prior achievement (*Math* $_0$); and (3) adding student characteristics and survey year. Results are shown in Figure D1.

**Figure D1**

*Standardized Coefficient Estimates of Sensemaking Constructs across Model Specifications* (± 1 SE) (*n = 3,170)*

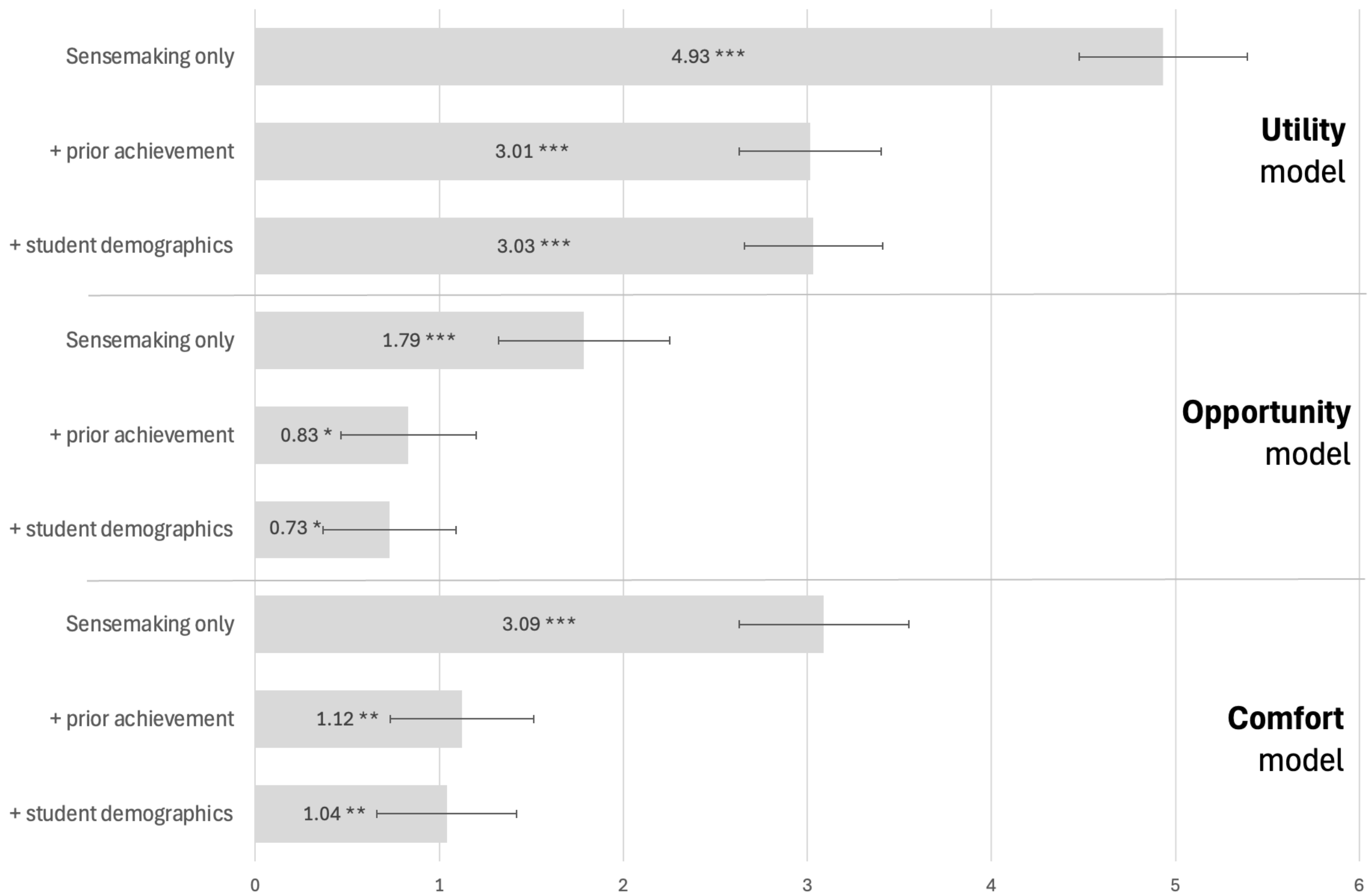

All three constructs retained significance across all specifications, demonstrating robust predictive validity. *Utility* showed the largest and most stable coefficients across specifications (*b* = 4.93, 3.01, 3.03), with modest attenuation from the unrestricted to the prior achievement model and virtually no further change with the addition of student demographics. *Comfort* showed a similar pattern of attenuation (*b* = 3.09 to 1.12 to 1.04), with a more pronounced drop once prior achievement was controlled, suggesting that a portion of unadjusted associations with student growth reflects prior achievement differences rather than independent effects. *Opportunity* showed the smallest coefficients across all specifications (b = 1.79, 0.83, 0.73) and retained only modest significance in the fully specified model ($p < .05$).

**Table D1**

*Predictive Validity of Sensemaking Constructs vs. Component Items (n = 3,170)*

| | Utility | | Opportunity | | Comfort | |
|---|---|---|---|---|---|---|
| | b | β | b | β | b | β |
| 1 **Solve** independently | 4.43 *** | 3.34 *** | -0.09 | 0.19 | 4.53 *** | 2.16 *** |
| 2 **Make** multiple representations | 2.36 *** | 1.49 *** | 2.18 * | 0.88 * | 1.99 * | 0.97 * |
| 3 **Compare** multiple representations | 1.19 * | 0.86 * | 0.27 | 0.34 | 0.78 | 0.67 |
| 4 **Explain** my strategy | 2.97 *** | 2.32 *** | 1.86 * | 0.89 * | 2.97 ** | 1.56 *** |
| 8 **Share** with whole class | 2.38 *** | 2.09 *** | 1.62 | 0.98 * | -0.06 | 0.54 |
| ***Sensemaking*** latent variable | | 3.03 *** | | 0.73 * | | 1.04 ** |

*Note.* Each listed item or construct was added individually to the fully specified model. * $p < .05$, ** $p < .01$, *** $p < .001$.
b = unstandardized coefficient; β = standardized coefficient